\documentclass[pdflatex,sn-mathphys-num]{sn-jnl}

\usepackage{graphicx} 
\usepackage{CJKutf8}
\usepackage[ruled,linesnumbered,vlined]{algorithm2e}
\usepackage{float} 
\usepackage{amscd}
\usepackage{amsmath,amssymb,amsfonts}%
\usepackage{mathtools}

\usepackage{amsthm}
\usepackage{threeparttable}
\usepackage{longtable}
\usepackage{caption}
\usepackage{subcaption}
\usepackage{booktabs}
\usepackage{adjustbox}
\usepackage{multirow}
\usepackage[numbers]{natbib}
\usepackage{setspace}
\usepackage{makecell}
\usepackage{bm}
\usepackage{calc} 
\usepackage{graphicx}%
\usepackage{mathrsfs}%
\usepackage[title]{appendix}%
\usepackage[utf8]{inputenc}

\newcommand{\vcenteredmultirow}[3][0pt]{\multirow{#2}{*}{\raisebox{-\height+\baselineskip+\heightof{#3}/2-\totalheightof{#3}/2+\baselineskip/2+#1}{#3}}}
\newcommand{\headerVcenter}[3][0pt]{
  \vcenteredmultirow[#1]{#2}{%
    \raisebox{0pt}[\heightof{#3}]{#3}
  }
}
\usepackage[marginal]{footmisc}

\renewcommand{\proofname}{\bf Proof}

\renewcommand{\figurename}{Fig.}
\SetKwInput{KwInput}{Input} 
\SetKwInput{KwOutput}{Output} 
\newcommand{\Step}[1]{\textbf{Step #1:}}

\newtheoremstyle{theorem}{3pt}{3pt}{\itshape}{0cm}{}{}{0.7em}{}
\theoremstyle{theorem}
\newtheorem{theorem}{\bf Theorem}
\newtheorem{lemma}{\bf Lemma}

\newtheoremstyle{remark}{3pt}{3pt}{}{0cm}{}{}{0.7em}{}
\theoremstyle{remark}
\newtheorem{remark}{\bf Remark}
\newtheorem{definition}{\bf Definition}%
\newtheorem{asp}{\bf Assumption}

\begin{document}
\captionsetup[figure]{labelfont={bf},name={Fig.},labelsep=space}

\title[Article Title]{Depth-first directional search for nonconvex optimization}

\author*[1]{\fnm{Yu-Xuan} \sur{Zhang}} \email{yx-z20@mails.tsinghua.edu.cn \thanks{12}}

\author[1]{\fnm{Wen-Xun} \sur{Xing}}\email{wxing@mail.tsinghua.edu.cn}

 \affil[1]{\orgdiv{Department of Mathematical Sciences}, \orgname{Tsinghua University}, \orgaddress{\city{Beijing}, \postcode{100084},  \country{China}}}

\abstract{
Random search methods are widely used for global optimization due to their theoretical generality and implementation simplicity. This paper proposes a depth-first directional search (DFDS) algorithm for globally solving nonconvex optimization problems. Motivated by the penetrating beam of a searchlight, DFDS performs a complete stepping line search along each sampled direction before proceeding to the next, contrasting with existing directional search methods that prioritize broad exploratory coverage. We establish the convergence and computational complexity of DFDS through a novel geometric framework that models the success probability of finding a global optimizer as the surface area
of a spherical cap. Numerical experiments on benchmark problems demonstrate that DFDS achieves significantly higher accuracy in locating the global optimum compared to other random search methods under the same function evaluation budget.}



\maketitle

\footnote{\noindent
        \textbf{Funding:} This work is supported by the National Natural Science Foundation of China Grant No. 11771243}

\section{Introduction}
\label{section introduction}

In this paper, we consider the following nonconvex optimization problem, 
\begin{equation}
\label{problem}
    \begin{array}{rl}
    &\min\quad  f(x)\\
    &{\rm s.t.} \quad x\in\mathcal{X},
    \end{array} \vspace{-1ex}
\end{equation}
where $f(x)$ is a real-valued continuous and nonconvex function, and $\mathcal{X}$ is a compact and convex set in the $N$-dimensional Euclidean space $\mathbb{R}^N$ with its diameter denoted by $D_0$. By the Weierstrass theorem, the global minimum value $f^*$ of (\ref{problem}) is attainable at some $x^* \in \mathcal{X}$, i.e., $f(x^*) = f^*$. Throughout this paper, we consider the following assumption: given an error threshold $\epsilon>0$, we can obtain a $R_\epsilon>0$ such that $f(x)\leq f(x^*)+\epsilon$ whenever $\|x-x^*\|\leq R_\epsilon$.

The primary focus of our research lies in globally solving optimization problems of the form (\ref{problem}) using random search methods. Random search constitutes a foundational class of algorithms in derivative-free optimization (DFO). DFO addresses problems where derivative information is unavailable, unreliable, or
impractical to obtain, such as legacy code integration, black-box simulations, or noisy function evaluations \cite{Rios2013,Conn2009,Cao2025}. 
Due to its applicability to non-differentiable systems and complex modeling scenarios, DFO has garnered significant attention in scientific computing, engineering and artificial intelligence \cite{Larson2019}. In the field of global optimization, derivative-free methods are particularly vital, as derivative-based approaches typically guarantee convergence only to stationary points \cite{Ghadimi2016,Cartis2010}, rather than the global optima.

Among DFO methods, random search is characterized by theoretical extensiveness and implementation simplicity. These methods typically impose minimal assumptions on the objective function $f$ and many exhibit global optimal properties by systematically exploring the entire feasible region. Based on their mechanisms of generating iterative points, we categorize random search methods into three classes, sampling methods, directional search methods and heuristic methods.

In typical random sampling methods, random points are sampled from a specific distribution, and the best solution is retained. Among these methods, pure random search (PRS) is the simplest: it sequentially generates independent and uniformly distributed points over $\mathcal{X}$, and terminates after a preset number of function evaluations. By sampling from the entire region, PRS is guaranteed to converge to the global minima with probability one. Furthermore, if $\mathcal{X}$ is an $N$-dimensional unit ball and $f$ satisfies the Lipschitz condition with constant $K$, the expected number of PRS sample points to achieve an $\epsilon-$optimal solution is $O((K/\epsilon)^N)$, which is exponential in $N$ \cite{Zabinsky2003_PRS}. Although PRS guarantees convergence theoretically, it suffers from the curse of dimensionality and fails to locate any global minima even in moderate dimensional problems \cite{Pepelyshev2018}.

To address PRS's lack of memory and independent sampling, adaptive search methods are developed. Pure adaptive search (PAS) generates a sequence of points uniformly distributed in the improving level set associated with each previous point. For a $f$ with its Lipschitz constant $K$ and feasible domain with diameter $D_0$, PAS achieves an $\epsilon-$optimal solution in $O(1+N[\ln(KD_0 / \epsilon)])$ expected improving iterations, exhibiting linear complexity in $N$ \cite{Zabinsky1992}. Other adaptive methods, including the hesitant adaptive search (HAS) \cite{Bulger1998} and annealing adaptive search (AAS) \cite{Romeijn1994_2}, are variants and relaxations of PAS. While adaptive search methods theoretically achieve linear complexity in convergence, their implementations are often infeasible as they rely on sampling distributions on improving level sets.

Directional random search methods iteratively sample random directions and generate iteration points via line search along these directions. Variants within this class primarily involve two aspects, either choosing the sampling distribution of directions or determining the step size selection strategy. Some methods enforce that sampled directions form a positive spanning set (PSS), whose conical combination spans $\mathbb{R}^N$. By coupling PSS with sufficiently small step sizes, descent conditions can be guaranteed. Recent theoretical advances have explored the convergence rates and worst-case complexity bounds of these algorithms to merely stationary points when applied to nonconvex $f$ \cite{Bergou2020,Bibi2020,Gratton2015}.

In contrast to PSS-based methods, some directional search methods sample random directions uniformly from a hypersphere. Among these methods, various step size selection strategies are suggested. For instance, a fixed step size method is proposed in \cite{Rastrigin1963} and later extended by \cite{Lawrence1972} to support bidirectional updates, and an adaptive step size method is introduced in \cite{Schumer1968}. As classified by \cite{Solis1981}, these methods constitute local search methods due to their bounded sampling support, which inherently precludes global optimality guarantees.

Hit-and-run methods \cite{Zabinsky2020} represent a specialized class of directional search methods designed for global optimization, practically realizing adaptive search methods via the hit-and-run generator. The hit-and-run sequence is generated by iteratively taking steps of random sizes along random directions and is guaranteed to asymptotically approach a uniform distribution over arbitrary sets \cite{Smith1984}. Motivated by PAS, the improving hit-and-run (IHR) generates a hit-and-run sequence of length 1 at each step and updates if the iteration point lands in the improving level set. For a class of elliptical programs, IHR achieves $\epsilon-$optimality in $O(N^{5/2})$ expected function evaluations \cite{Zabinsky1993}. However, this complexity result has limited applicability, as it is derived for a highly specific case where the function value depends solely on a matrix-induced norm measuring the distance to the optimal point. As an variant of IHR, Hide-and-seek incorporates an acceptance probability with a cooling schedule, enabling convergence to the global minima in general problems, though its computational complexity remains unquantified \cite{Romeijn1994_1}.

Heuristic methods explore the feasible region through diverse generation and update mechanisms. For instance, the Metropolis rule is adopted in simulated annealing (SA) and the pseudo-random proportional rule governs decision-making in ant colony system (ACS). Heuristic algorithms are relatively straightforward in implementation, and some have been proven to converge to the global minima in probability 1 under specific conditions \cite{Locatelli2002,Xu2018}. However, systematic analyses of their computational complexity remain limited in the literature.

We conclude that within random search methods, sampling methods possess extensive theoretical analyses, however, their implementations often suffer from inefficiency. In contrast, heuristics and direction search methods are designed for practical considerations, with some proven to converge to the global minima in probability 1, yet complexity analyses for these methods remain scarce.

In this paper, we propose a random search method which systematically explores the feasible domain by generating random directions, and performing a uniform stepping line search along each direction until reaching the boundary. Our full-depth line search design is inspired by the unidirectional penetration of a searchlight beam. Unlike existing random search methods that prioritize broad coverage (either sampling points directly or evaluating a single point along each direction), our approach adopts a depth-first strategy, exhaustively traversing each random direction before switching to the next. We term this method the depth-first directional search (DFDS).

We investigate on the global $\epsilon-$optimal properties of the DFDS algorithm for a given error threshold $\epsilon$. As previously stated, it is assumed that the set of $\epsilon-$optimal points contains a ball centered at the global minimizer. Under this assumption, the set of directions enabling DFDS to locate an $\epsilon-$optimal solution via line search corresponds to a spherical cap on the unit sphere. We derive a lower bound on the relative surface area of this spherical cap and, consequently establish the convergence in probability and computational complexity of DFDS. In particular, we show that DFDS attains an $\epsilon-$accuracy with at most $O((\frac{2(D_0+R_\epsilon)}{\sqrt{3}R_\epsilon})^{N}\,\frac{D_0}{R_\epsilon}\,\sqrt{N}\,$\scalebox{1.25}{$\frac{1}{\epsilon}$}$)$ function evaluations required in expectation. By introducing a geometric framework, our theoretical analysis provide a novel perspective for establishing convergence and complexity bounds in directional search methods. 

Finally, we conduct numerical experiments to evaluate the performance of DFDS on benchmark nonconvex optimization problems, comparing it against two random search methods, PRS and IHR. The selection of baselines is motivated by two key considerations: PRS and DFDS both belong to the exponential complexity class; PRS and IHR are both directional search methods while they diverge in their line search strategies. Experimental results show that DFDS achieves not only higher accuracy than PRS and IHR under the same evaluation budget in most instances, but also exhibits superior scalability, with its accuracy advantage becoming more pronounced as the problem dimension increases.

The rest of the paper is organized as follows. In Section 2, notations, assumptions and some relevant results are provided. In Section 3, we present DFDS and state theoretically the convergence and computational complexity of DFDS. Numerical tests are discussed in Section 4. Finally, concluding remarks are presented in Section 5.

\section{Preliminaries}
\label{section preliminaries}

General notations are listed as follows. $\mathbb{R}^N$ denotes the $N-$dimensional Euclidean space and $\|\cdot\|$ denotes the Euclidean norm. For points $x,y$ in $\mathbb{R}^N$, $d(x,y)=\|x-y\|$ denotes their Euclidean distance. For vectors $x,y$ in $\mathbb{R}^N$, $\langle x,y\rangle$ denotes their included angle and $(x, y)$ denotes their inner product. For a point $x$ and a vector $d$ both in $\mathbb{R}^N$, $\mathcal{R}(x,d)=\{x+td | t\geq0\}$ denotes the ray starting from $x$ along $d$. A mathematical calligraphy font letter (e.g., $\mathcal{S}$) denotes a set. For a set $\mathcal{S}$, $d_{\mathcal{S}}(x)=\min\limits_{s\in\mathcal{S}}d(x,s)$ denotes the Euclidean distance from a point $x$ to $\mathcal{S}$ and $d(\mathcal{S})=\max\limits_{s_1,s_2\in\mathcal{S}} \|s_1-s_2\|$ denotes the diameter of $\mathcal{S}$. $m(\mathcal{S})$ denotes the measure of $\mathcal{S}$. $\mathcal{B}(x,R)=\{y\in \mathbb{R}^N|d(x,y)<R\}$ denotes the open ball centered at $x$ with radius $R$ and $\overline{\mathcal{B}}(x,R)$ denotes its closure; $\mathcal{S}^{N-1}=\{x\in\mathbb{R}^N|\|x\|= 1\}$ denotes the ($N$-1)-dimensional sphere.

Furthermore, some definitions and relevant results are given below.

\begin{definition}
    \label{extend set}
    \rm{(Extended Set.) Given a set $\mathcal{X}$ in $\mathbb{R}^N$ and a $R>0$, the $R$-extended set $\mathcal{X}_R$ of $\mathcal{X}$ is defined as 
    \begin{equation*}
        \mathcal{X}_R=\{x\in\mathbb{R}^N|d_{\mathcal{X}}(x)\leq R\}.
    \end{equation*} }
\end{definition}

\begin{definition}
    \label{epsilon better}
    \rm{($\epsilon-$better and $\epsilon-$optimal.) Given an error threshold $\epsilon>0$ and any $x_1\in\mathcal{X}$, $x_2\in\mathcal{X}$ is called an $\epsilon-$better point than $x_1$ if $f(x_2)\leq f(x_1)-\epsilon$; $x$ is called an $\epsilon-$optimal point if $f(x)\leq f^*+\epsilon$, where $f^*$ is the global minimum value of (\ref{problem}).}
\end{definition}

Finally, the following assumptions are considered throughout the paper. 

\begin{asp}
    \label{asp XR}
    \rm{Given any $R>0$, it is trivial to verify whether $x$ is in $\mathcal{X}_R$. }
\end{asp}
\begin{asp}
\label{asp R0}
    \rm{For a given threshold $\epsilon>0$, there exists a $R_\epsilon>0$ such that $f$ is well-defined on $\mathcal{X}_{R_\epsilon}$ and for every $\epsilon-$optimal point $x\in\mathcal{X}_{R_\epsilon}$,
    \begin{center}
        $f(y) \leq f(x) + \epsilon/3$
    \end{center}
    holds for all  $ y\in\overline{\mathcal{B}}(x,R_\epsilon) \cap\mathcal{X}_{R_\epsilon}$.}
\end{asp}
The existence of $R_\epsilon$ in Assumption \ref{asp R0} is guaranteed for any continuous $f$ on a compact $\mathcal{X}$. By the Heine–Cantor theorem, $f$ is uniformly continuous on any $\mathcal{X}_{R}$ (which is also compact). Hence, for any $\epsilon>0$, there exists a $\delta>0$ such that $|f(x)-f(y)|\leq\frac{\epsilon}{3}$ holds for all $x,y\in\mathcal{X}_R$ with $\|x-y\|\leq\delta$. Taking $R_\epsilon=\min\{R,\delta\}$ therefore satisfies Assumption \ref{asp R0}. Further discussions on the derivation of $R_\epsilon$ are provided at the end of Section \ref{section main}.

For some general $\mathcal{X}$ and $f$, the global minimum may be attained on the boundary of $\mathcal{X}$. The algorithm is designed on $\mathcal{X}_{R_\epsilon}$ to include the $R_\epsilon-$neighborhood of each global minimum.

\section{Depth-First Directional Search}
\label{section main}

In this section, we present the depth-first directional search (DFDS) algorithm and discuss the convergence and computational complexity of DFDS. The algorithm is defined as follows. \vspace{2ex}

\begin{algorithm}[H]
\caption{\textbf{Depth-First Directional Search (DFDS)}}
\KwInput{ An initial step size $R_\epsilon$, the number of generated directions $M$ and an error tolerance $\epsilon$.}
\KwOutput{A feasible point of (\ref{problem}).}

\Step{0} Select $x_0\in\mathcal{X}$ and set $k=0$, $m=0$.

\Step{1} \Indp If $m\geq M$, go to \textbf{Step 3}. \\
 Otherwise generate a $d_m$ from the uniform distribution on $\mathcal{S}^{N-1}$ and set $r=R_\epsilon$. \Indm

\Step{2} If $x_k+rd_m\notin\mathcal{X}_{R_\epsilon}$, set $m=m+1$ and go to \textbf{Step 1}.\\
Else if $f(x_k+rd_m)\leq f(x_k)-\epsilon/3$, set $x_{k+1} = x_k+rd_m$, set $k=k+1$, set $m=0$ and go to \textbf{Step 1}.\\
Otherwise let $r= r+R_\epsilon$ and go to \textbf{Step 2}.

\Step{3} Select a $\widetilde{x_k}$ in $\mathcal{X}\cap\overline{B}(x_k,R_\epsilon)$ and output $\widetilde{x_k}$.
\end{algorithm} \vspace{2ex}

The algorithm proceeds as follows. At each iteration point $x_k$ $(k\geq0)$, it generates at most $M$ random directions on $\mathcal{S}^N$. Along each direction $d_m$, it performs multiple \textbf{function evaluations} (defined as any point where the objective function is computed) at a uniform spacing of $R_\epsilon$ within the extended set $\mathcal{X}_{R_\epsilon}$. The \textbf{iteration point} $x_k$ is updated if and only if an $\epsilon/3-$better evaluation point is obtained. The entire process, starting from the initial point $x_0$ and continuing until no $\epsilon/3-$better point is found after sampling $M$ directions at some iteration point $x_k$, is termed the \textbf{directional random search phase}. After the directional random search phase terminates, the output point $\widetilde{x_k}$ is selected from $\mathcal{X}\cap\overline{\mathcal{B}}(x_k,R_\epsilon)$ as a feasible candidate of $x_k$. The attainability of $\widetilde{x_k}$ is ensured by Assumption \ref{asp XR}. \vspace{1ex}

 The theoretical analysis in this section begins by addressing the $\epsilon-$optimality of DFDS (i.e., the output point $\widetilde{x_k}$ being $\epsilon-$optimal). A sufficient condition for this $\epsilon-$optimality is given by the following lemma.
\begin{lemma}
    \label{lemma optimality}
    For any initial point $x_0\in\mathcal{X}$, DFDS terminates within $k_{max}\coloneqq \big\lfloor \frac{f(x_0)-f^*}{\epsilon/3}\big\rfloor$ iteration points. Moreover, DFDS outputs an $\epsilon-$optimal point in $\mathcal{X}$ if $x_k$ is $\frac{2\epsilon}{3}-$optimal for some $k\geq0$.
\end{lemma}

\begin{proof}
 In DFDS, each iteration reduces the function value by at least $\epsilon/3$. Hence, the number of iteration points is bounded by $k_{max}\coloneqq \big\lfloor \frac{f(x_0)-f^*}{\epsilon/3}\big\rfloor$ given any initial point $x_0\in\mathcal{X}$.
 
Denote by $K$ the index of the iteration point at which the directional random search phase terminates, then $K\leq \big\lfloor\frac{f(x_0)-f^*}{\epsilon/3}\big\rfloor$.
If $ \exists\, k\geq 0 $ such that $x_k$ is $\frac{2\epsilon}{3}-$optimal, then \vspace{-1ex}
 \begin{equation*}
        f(\widetilde{x_K})\leq f(x_K)+\frac{\epsilon}{3} \leq f(x_k)+\frac{\epsilon}{3} \leq f^*+\epsilon  
\end{equation*}  
  according to Assumption \ref{asp R0} and the $\epsilon-$optimality of DFDS is ensured.
\end{proof} 

For notation convenience, we introduce the following definition.
\begin{definition}
     \rm{For the given threshold $\epsilon>0$, $\mathcal{X}_{R_\epsilon}^*\coloneqq\{x\in\mathcal{X}_{R_\epsilon}|f(x)\leq f^*+\frac{2\epsilon}{3}\}$ is defined as the $\frac{2\epsilon}{3}-$optimal set in $\mathcal{X}_{R_\epsilon}$.}
\end{definition} 

Following Lemma \ref{lemma optimality}, we analyze the probability that the directional random search phase reaches an iteration point in $\mathcal{X}^*_{R_\epsilon}$. Such analysis is started by considering the directional search at some fixed $x_k\notin\mathcal{X}^*_{R_\epsilon}$. At $x_k\notin\mathcal{X}^*_{R_\epsilon}$, the process of uniformly sampling a $d_m$ and conducting line search along $d_m$ is referred to as one random trial. In each random trial at $x_k$, the probability of finding an $\frac{\epsilon}{3}-$better point than $x_k$ is identical.

\begin{definition}
    At $x_k$, $p_{x_k}$ is defined as the probability of finding an $\frac{\epsilon}{3}-$better point than $x_k$ in a single random trial. This probability remains identical across all trials at $x_k$.
\end{definition}

 A lower bound for $p_{x_k}$ will be established in Theorem \ref{theo p}. To this end, we first present the following necessary definitions and preliminary lemmas.

\begin{lemma}
    \label{lemm set}
    Let $x^*\in\mathcal{X}$ be a global minimizer of (\ref{problem}) and $x_k\notin\mathcal{X}^*_{R_\epsilon}$. Every point in $\overline{\mathcal{B}}(x^*,R_\epsilon)$ is $\frac{\epsilon}{3}-$optimal, and thus $\frac{\epsilon}{3}-$better than $x_k$.
\end{lemma}

\begin{proof}
    According to Assumption \ref{asp R0}, $f(x)\leq f(x^*)+\frac{\epsilon}{3}$ for all $ x\in\overline{\mathcal{B}}(x^*,R_\epsilon)$. Since $x_k\notin\mathcal{X}^{*}_{R_\epsilon}$, we have \vspace{-1ex}
    \begin{equation*}
        f(x)\leq f(x^*)+\frac{\epsilon}{3}=f^*+\frac{\epsilon}{3}< f(x_k)- \frac{\epsilon}{3} \vspace{-1ex}
    \end{equation*} 
    for all $ x\in\overline{\mathcal{B}}(x^*,R_\epsilon)$, which completes the proof.
\end{proof}

\begin{definition}
Let $x^*\in\mathcal{X}$ be a global minimizer of (\ref{problem}) and $x_k\notin\mathcal{X}^*_{R_\epsilon}$. The angle $\alpha$ is defined as $\alpha=\arcsin\frac{\sqrt{3}R_\epsilon}{2||x^*-x_k||}$ and the set of directions is defined as $\mathcal{S}_{\alpha}^{N-1}=\{d\in\mathcal{S}^{N-1}| \langle x^*-x_k,d\rangle \leq\alpha\}$.
\end{definition}

By Lemma \ref{lemm set}, $x_k\notin\overline{\mathcal{B}}(x^*,R_\epsilon)$ and thus $\alpha$ is well-defined. In the definition above, $\alpha$ and $\mathcal{S}_\alpha^{N-1}$ are used for notational convenience, omitting their dependence on $x_k$ and $x^*$. Specifically, $\mathcal{S}_\alpha^{N-1}$ represents a spherical cap with the colatitude angle $\alpha$. The following lemma characterizes key properties of $\alpha$ and $\mathcal{S}_\alpha^{N-1}$.

\begin{lemma}
\label{lemm alpha and D}
Let $x^*\in\mathcal{X}$ be a global minimizer of (\ref{problem}) and $x_k\notin\mathcal{X}^*_{R_\epsilon}$. The angle $\alpha$ and the set $\mathcal{S}_\alpha^{N-1}$ satisfies: \\ \vspace{1.2ex}
 \indent 1. $\frac{\sqrt{3}R_\epsilon}{2(D_0+R_\epsilon)}\leq\sin\alpha<\frac{\sqrt{3}}{2}$, and equivalently $\arcsin\big(\frac{\sqrt{3}R_\epsilon}{2(D_0+R_\epsilon)}\big)\leq\alpha<\frac{\pi}{3}$. \vspace{1.2ex} \\
 \indent 2. For every $d\in\mathcal{S}_\alpha^{N-1}$, ray $\mathcal{R}(x_k,d)$ intersects $\overline{\mathcal{B}}(x^*,R_\epsilon)$, forming a chord of length at least $R_\epsilon$. \\
 \indent 3. If any generated $d_m$ is in $\mathcal{S}_\alpha^{N-1}$, then the directional search at $x_k$ succeeds to find an $\frac{\epsilon}{3}-$better point for certain.
\end{lemma}
\begin{proof}
    Since $x_k\in\mathcal{X}_{R_\epsilon}$ and $x^*\in\mathcal{X}$, $\|x^*-x_k\|\leq d(\mathcal{X})=D_0+R_\epsilon$, while $\|x^*-x_k\|>R_\epsilon$ is given in Lemma \ref{lemm set}. Together with $\sin\alpha = \frac{\sqrt{3}R_\epsilon}{2\|x^*-x_k\|}$, the first conclusion is obtained.
    
    The distance from $x^*$ to ray $\mathcal{R}(x_k,d)$ equals $\|x^*-x_k\|\cdot\sin\big(\langle x^*-x_k,d\rangle\big)$, and is no greater than $\frac{\sqrt{3}}{2}R_\epsilon$ for every $d\in\mathcal{S}_\alpha^{N-1}$. Therefore 
    $\mathcal{R}(x_k,d)$ intersects $\overline{B}(x^*,R_\epsilon)$ with a chord of length at least $R_\epsilon$, which is illustrated in Fig. \ref{fig lineseg}. 
    \begin{figure}
        \centering
        \includegraphics[width=.65\textwidth]{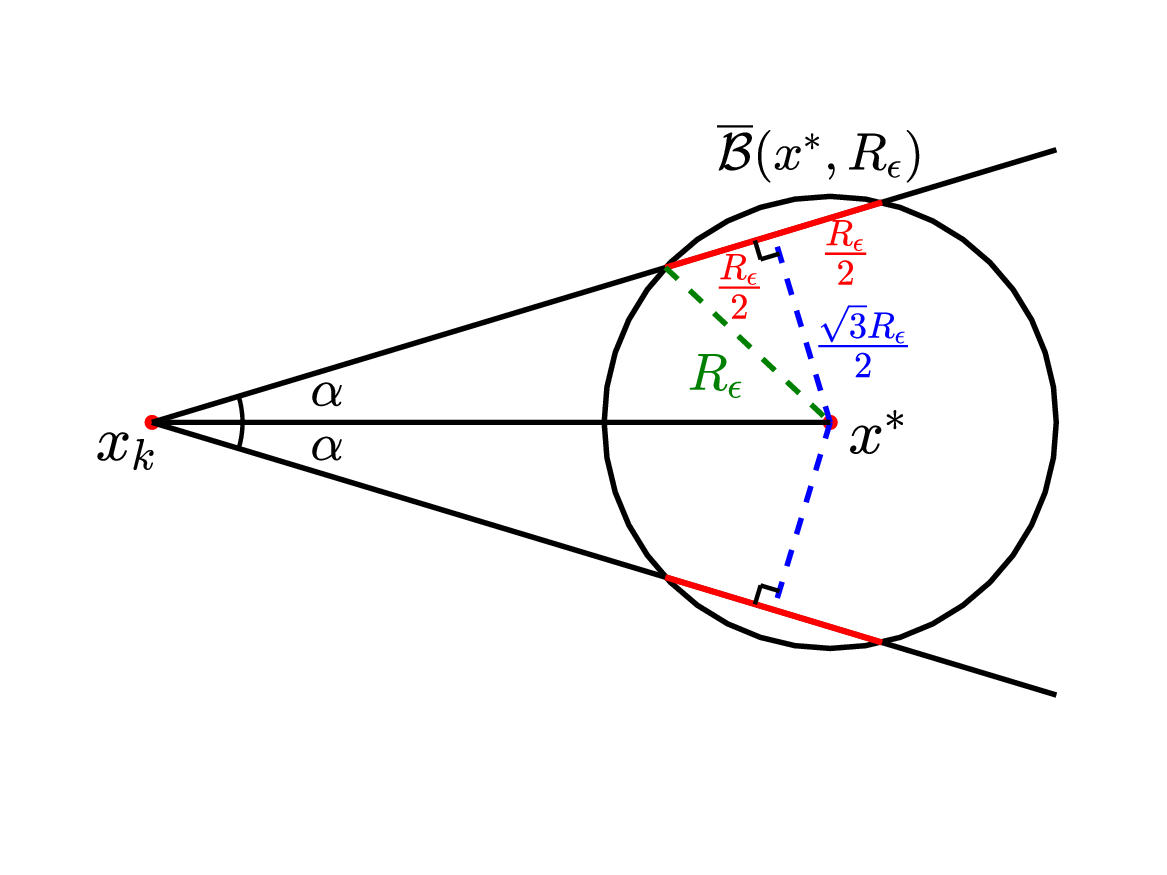}
        \caption{$R(x_k,d)$ when $\langle x^*-x_k,d\rangle=\alpha$ and $\overline{B}(x^*,R_\epsilon)$ ($N=2$)}
        \label{fig lineseg}
    \end{figure} 

    The algorithm uniformly searches points at a distance of $R_\epsilon$ along each direction $d_m$ starting from $x_k$. For any $d_m\in\mathcal{S}_\alpha^{N-1}$, the intersection chord $\mathcal{R}(x_k,d)\cap\overline{\mathcal{B}}(x^*,R_\epsilon)$ guarantees an $\frac{\epsilon}{3}-$better point obtained along $d_m$.
\end{proof}

Following the analysis above, it can be concluded that the probability of a generated $d_m$ belonging to $\mathcal{S}_{\alpha}^{N-1}$ serves as a lower bound for $p_{x_k}$. Such a probability is estimated in the following theorem.

\begin{theorem}
\label{theo p}
Let $x^*\in\mathcal{X}$ be a global minimizer of (\ref{problem}) and $x_k\notin\mathcal{X}^*_{R_\epsilon}$. For a $d_m$ generated from the uniform distribution on $\mathcal{S}^{N-1}$, define $p_{N,\alpha}$ as the probability that $d_m\in\mathcal{S}_{\alpha}^{N-1}$. The following statements hold for $p_{N,\alpha}$, \\
    1. \begin{equation*}
\displaystyle p_{N,\alpha}=
\left\{
    \begin{array}{lr}
        \displaystyle\frac{1}{\pi} \, \alpha\,  &  N=2 \vspace*{1.2ex} \\
        \displaystyle\frac{1}{2}\big(1-\cos\alpha\big)\,  & N=3  \vspace*{1.2ex} \\      
        \displaystyle\frac{1}{\pi}\Big(\alpha - \sum\limits_{t=0}^{\frac{N}{2}-2} \frac{(2t)!!}{(2t+1)!!}\cos\alpha\sin^{2t+1}\alpha\Big)  & \text{N is even and }N\geq4 \vspace*{1.2ex}  \\   
        \displaystyle\frac{1}{2} \cdot \Big[ 1-\cos\alpha\Big(1+\sum_{t=1}^{\frac{N-3}{2}} \frac{(2t-1)!!}{(2t)!!}\sin^{2t}\alpha\Big)\Big] & \text{N is odd and }N\geq5.
    \end{array}
\right.
\end{equation*}
2. \begin{equation*}
    \displaystyle p_{N,\alpha}\geq
    \left\{
    \begin{array}{ll}
    \displaystyle\frac{(N-2)!!}{\pi N\cdot (N-3)!!} \,\cos\alpha\sin^N\alpha  & \text{N is even and }N\geq4  \vspace*{1.2ex} \\
        \displaystyle\frac{(N-2)!!}{2(N-1)!!} \,\cos\alpha\sin^{N-1}\alpha & \text{N is odd and }N\geq5. 
    \end{array}
    \right.
\end{equation*}
\end{theorem}
\begin{proof}
\begin{figure}
    \centering
    \includegraphics[width=.55\textwidth]{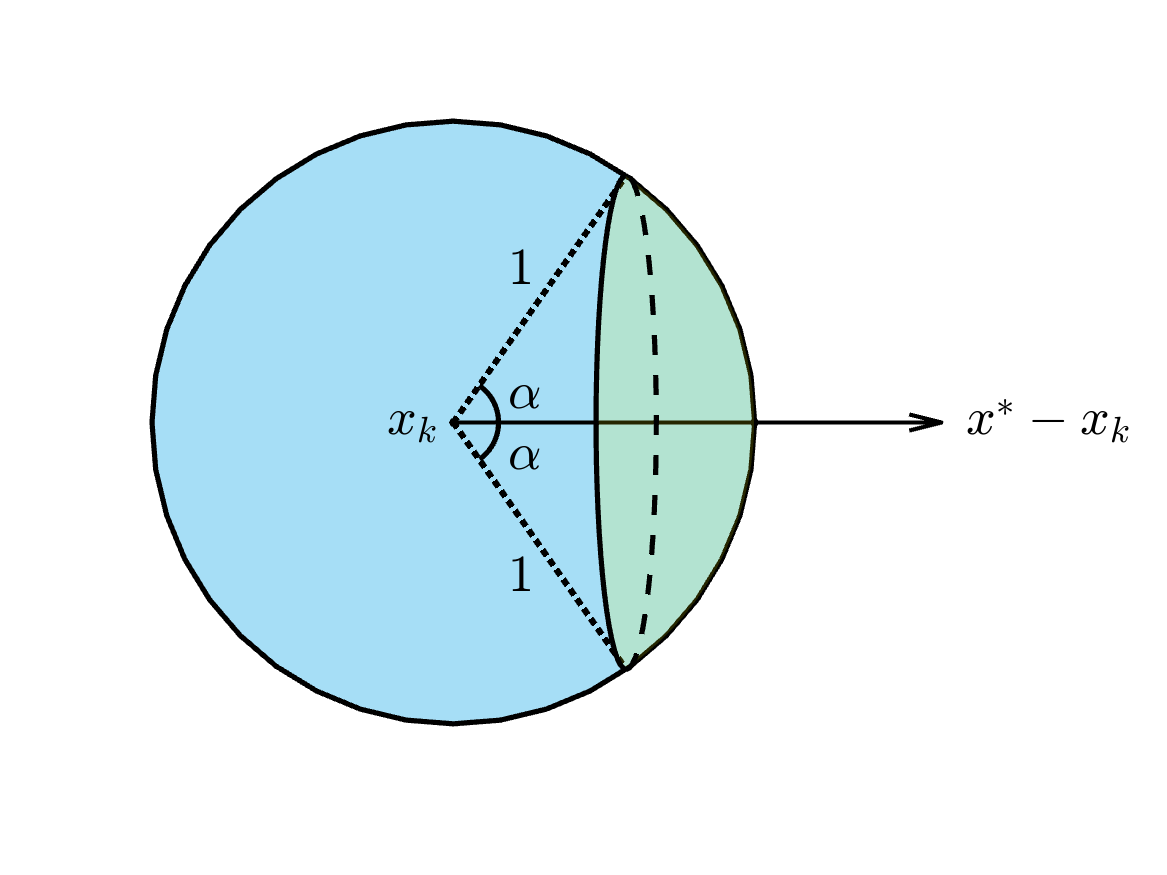}
    \caption{The spherical cap $\mathcal{S}_{\alpha}^2$ of the unit ball in $\mathbb{R}^3$}
    \label{fig cap}
    \end{figure}
\indent For a $d_m$ generated from the uniform distribution on the unit sphere $\mathcal{S}^{N-1}$, $p_{N,\alpha}$ can be computed as the proportion of the area of $\mathcal{S}_\alpha^{N-1}$ to that of $\mathcal{S}^{N-1}$, \vspace{-1ex}
\begin{equation}
    \label{geometric probability}
    p_{N,\alpha}=\frac{m(\mathcal{S}_\alpha^{N-1})}{w_{N-1}} \vspace{-1ex}
\end{equation}

\noindent where $m(\mathcal{S}_\alpha^{N-1})$ and $w_{N-1}$ denote the surface areas of $\mathcal{S}_\alpha^{N-1}$ and $\mathcal{S}^{N-1}$ respectively. For illustration, we demonstrate $\mathcal{S}_\alpha^2$ (the green shaded area) and $\mathcal{S}^2$ (the blue shaded area and the green shaded area as a whole) in Fig. \ref{fig cap}. The surface area formula of the spherical cap $\mathcal{S}_\alpha^{N-1}$ \cite{Li2011} is given by
\begin{equation*}
\label{cap area}
m( \mathcal{S}^{N-1}_\alpha) = w_{N-2}\int_{0}^{\alpha} \sin^{N-2}\theta \,d\theta .
\end{equation*}     
\noindent Therefore we have
\begin{equation}
\label{general probability}
    \begin{aligned}
        p_{N,\alpha} = \frac{w_{N-2}}{w_{N-1}}\int_0^{\alpha} \sin^{N-2}\theta \,d\theta  .
    \end{aligned}
\end{equation}
For any $N\geq 2$, $w_{N-1}$ is given by  \vspace{-1ex}
\begin{equation}
\label{sphere area}
w_{N-1}=\frac{2\pi^{N/2}}{\Gamma(N/2)},
\end{equation}
where $\Gamma(\cdot)$ is the Gamma function with 
\begin{equation}
\label{eq gamma}
\displaystyle\Gamma(\frac{N}{2})=
\left\{
\begin{array}{ll}
    \displaystyle (\frac{N}{2}-1)\,(\frac{N}{2}-2)\cdots1  \vspace*{1.2ex} & N\text{ is even } \\
    \displaystyle  (\frac{N}{2}-1)\,(\frac{N}{2}-2)\cdots\frac{1}{2}\sqrt{\pi} \quad &  N\text{ is odd }.
\end{array}
\right.
\end{equation}
\noindent When $N\geq2$,  \vspace{-1ex}
\begin{equation*}
    \begin{array}{lrll}
    & \displaystyle\int_0^{\alpha} \sin^{N}\theta d\theta   &= \ &-\cos\alpha\sin^{N-1}\alpha - \displaystyle\int_0^{\alpha} (N-1)(\sin^2\theta-1)\sin^{N-2}\theta \, d\theta \vspace{2ex} \\
    \Longrightarrow\; &
    N \displaystyle \int_0^{\alpha}  \sin^{N}\theta \, d\theta & = \ & -\cos\alpha\sin^{N-1}\alpha +(N-1) \displaystyle \int_0^{\alpha} \sin^{N-2}\theta \, d\theta . 
\end{array}     \vspace{-1ex}
\end{equation*}
Denote $I_{N,\alpha} = \displaystyle\int_0^{\alpha} \sin^{N}\theta d\theta$, we have when $N\geq2$  \vspace{-1ex}
\begin{equation}
    \label{integral by parts}
        { I_{N,\alpha}}  = \  { \frac{N-1}{N}  \ I_{N-2,\alpha} - \frac{1}{N} \  \cos\alpha\sin^{N-1} \alpha};      
\end{equation}
and when $N=0,1$,
\begin{equation}
\label{eq I01}
    \begin{aligned}
        I_{1,\alpha} = \ & \int_0^{\alpha} \sin\theta \,d\theta 
        =  1-\cos\alpha \quad\text{and}\quad
        I_{0,\alpha} = \  \int_0^{\alpha} \sin^0\theta \,d\theta 
        =  \alpha.
    \end{aligned}
\end{equation}
Considering that $\displaystyle I_{N,\alpha}$ and $\displaystyle\Gamma(N/2)$ have different expressions when $N$ is odd or even, the calculation of $p_{N,\alpha}$ is categorized by the parity of $N$.\\
\indent Condition 1: $N$ is even. \\
\indent When $N=2$, $\displaystyle p_{2,\alpha}=\frac{m(\mathcal{S}_\alpha^{1})}{w_1}=\frac{2\alpha}{2\pi}=\frac{1}{\pi}\alpha$ is directly calculated from (\ref{geometric probability}).\\
\indent  When $N\geq 4$, for the first part of (\ref{general probability}) we have from (\ref{sphere area}) and (\ref{eq gamma}),
\begin{equation}
    \label{even w}
    \begin{aligned}
        \frac{w_{N-2}}{w_{N-1}} 
        &= \frac{\Gamma(\frac{N}{2})}{\sqrt{\pi}\,\Gamma(\frac{N-1}{2})}   \\[1.3mm]
        & = \frac{ (\frac{N}{2}-1)\,(\frac{N}{2}-2)\cdots1}{\sqrt{\pi}\, (\frac{N-1}{2}-1)\,(\frac{N-1}{2}-2)\cdots\frac{1}{2}\,\sqrt{\pi}}  \\[1.3mm]
        & = \frac{1}{\pi}\,\frac{(\frac{1}{2})^{N/2-1}\,(N-2)\,(N-4)\cdots2}{(\frac{1}{2})^{N/2-1}\, (N-3) \, (N-5)\cdots1}   \\[1.3mm]
        & = \frac{(N-2)!!}{\pi (N-3)!!} \, . 
    \end{aligned}
\end{equation}

\noindent When $N\geq4$, for the second part of (\ref{general probability}), we have from (\ref{integral by parts}),
\begin{equation}
    \label{even int}
    \begin{aligned}
          I_{N-2,\alpha} \ & = \tfrac{N-3}{N-2}\ I_{N-4,\alpha} - \tfrac{1}{N-2} \cos\alpha\sin^{N-3}\alpha  \\[1.3mm] 
    & = \ \tfrac{(N-3)(N-5)}{(N-2)(N-4)}\ I_{N-6,\alpha} - \tfrac{N-3}{N-2}\tfrac{1}{N-4} \cos\alpha\sin^{N-5}\alpha - \tfrac{1}{N-2} \cos\alpha\sin^{N-3}\alpha  \\[1.3mm]
    & =   \cdots \\[1.3mm]
        &\; \begin{aligned} 
            = \ & \tfrac{(N-3)(N-5)\cdots 1}{(N-2)(N-4)\cdots 2} \ I_{0,\alpha} - \tfrac{(N-3)(N-5)\cdots 3}{(N-2)(N-4)\cdots 4} \tfrac{1}{2} \cos\alpha\sin\alpha - \cdots \\[1.3mm] 
            & -\tfrac{N-3}{N-2} \tfrac{1}{N-4}\cos\alpha\sin^{N-5}\alpha -  \tfrac{1}{N-2} \cos\alpha\sin^{N-3}\alpha
    \end{aligned}   \\[1.3mm]
    & \; \begin{aligned}
        = \ & \tfrac{(N-3)!!}{(N-2)!!}\ \alpha - \sum\limits_{t=1}^{N/2-2} \tfrac{(N-3)(N-5)\cdots(2t+1)}{(N-2)(N-4)\cdots(2t+2)}\tfrac{1}{2t}\cos\alpha\sin^{2t-1}\alpha \\[1.3mm]
     & - \tfrac{1}{N-2} \cos\alpha\sin^{N-3}\alpha. 
    \end{aligned}    
    \end{aligned}
\end{equation}

\noindent Combining (\ref{even w}) and (\ref{even int}), when $N\geq 4$ we have
\begin{equation}
    \label{even p}
    \begin{aligned}
        p_{N,\alpha} \ & = \ \tfrac{w_{N-2}}{w_{N-1}}\ I_{N-2,\alpha} \\[1mm]
        &\; 
        \begin{aligned} 
        = \  \tfrac{(N-2)!!}{\pi(N-3)!!}\Big[&\tfrac{(N-3)!!}{(N-2)!!}\ \alpha - \sum\limits_{t=1}^{N/2-2} \tfrac{(N-3)(N-5)\cdots(2t+1)}{(N-2)(N-4)\cdots(2t+2)}\tfrac{1}{2t}\cos\alpha\sin^{2t-1}\alpha  \\[1mm]
        & - \tfrac{1}{N-2} \cos\alpha\sin^{N-3}\alpha \Big]
        \end{aligned}
        \\[1mm]
        & \; \begin{aligned}
            =\frac{1}{\pi}\big[ &\alpha -\cos\alpha\sin\alpha -\tfrac{2}{3}\cos\alpha\sin^3\alpha-\tfrac{4\cdot2}{5\cdot3}\cos\alpha\sin^5\alpha-\cdots \\[1.3mm]
        &-\tfrac{(N-4)!!}{(N-3)!!}\cos\alpha\sin^{N-3}\alpha\big]
        \end{aligned}
        \\[1.3mm]
        & = \frac{1}{\pi}\big(\alpha - \sum\limits_{t=0}^{N/2-2} \tfrac{(2t)!!}{(2t+1)!!}\cos\alpha\sin^{2t+1}\alpha\big).
    \end{aligned}
\end{equation}

\indent Condition 2: $N$ is odd. \\[2mm]
\indent When $N=3$, $\displaystyle p_{N,\alpha}=\frac{w_1}{w_2}\, I_{1,\alpha}=\frac{1}{2}\big(1-\cos\alpha\big)$ is calculated from (\ref{general probability}) and (\ref{eq I01}). \\ [2mm]
\indent When $N\geq5$, for the first part of (\ref{general probability}) we have from (\ref{sphere area}) and (\ref{eq gamma}),
\begin{equation}
    \label{odd w}
        \begin{aligned}
            \frac{w_{N-2}}{w_{N-1}} \
        & = \frac{\Gamma(\frac{N}{2})}{\sqrt{\pi}\,\Gamma(\frac{N-1}{2})}   \\[1.3mm]
        & = \ \frac{ (\frac{N}{2}-1)\,(\frac{N}{2}-2)\cdots\frac{3}{2}\cdot\frac{1}{2} \, \sqrt{\pi}}{\sqrt{\pi}\, (\frac{N-1}{2}-1)\,(\frac{N-1}{2}-2)\cdots1} \\[1.3mm]
        & = \ \frac{(\frac{1}{2})^{(N-1)/2}\,(N-2)\,(N-4)\cdots3\cdot1}{(\frac{1}{2})^{(N-1)/2-1}\, (N-3) \, (N-5)\cdots2}  \\[1.3mm]
        & = \ \frac{(N-2)!!}{2 (N-3)!!} \, .
        \end{aligned}
    \end{equation}

\noindent When $N\geq5$, for the second part of (\ref{general probability}), we have from (\ref{integral by parts}),
\begin{equation}
    \label{odd int}
    \begin{aligned}
        I_{N-2,\alpha} \ & = \ \tfrac{N-3}{N-2}\ I_{N-4,\alpha} - \tfrac{1}{N-2} \cos\alpha\sin^{N-3}\alpha \\[1.3mm]
        & = \  \tfrac{(N-3)(N-5)}{(N-2)(N-4)}\ I_{N-6,\alpha} - \tfrac{N-3}{N-2}\tfrac{1}{N-4}\cos\alpha\sin^{N-5}\alpha  - \tfrac{1}{N-2} \cos\alpha\sin^{N-3}\alpha  \\[1.3mm]
        & = \  \cdots \\[1.3mm]
        & \;
        \begin{aligned} 
            = \ &\tfrac{(N-3)(N-5)\cdots2}{(N-2)(N-4)\cdots 3}\ I_{1,\alpha}- \tfrac{(N-3)(N-5)\cdots4}{(N-2)(N-4)\cdots 5}\tfrac{1}{3}\cos\alpha\sin^2\alpha- \cdots \\[1.3mm]
                &- \tfrac{N-3}{(N-2)}\tfrac{1}{N-4}\cos\alpha\sin^{N-5}\alpha - \tfrac{1}{N-2} \cos\alpha\sin^{N-3}\alpha   
        \end{aligned}   \\[1.3mm]
        & \begin{aligned}
           = \ & \tfrac{(N-3)!!}{(N-2)!!}\,(1-\cos\alpha )-\sum\limits_{t=1}^{(N-5)/2}\tfrac{(N-3)(N-5)\cdots(2t+2)}{(N-2)(N- 4)\cdots(2t+3)}\tfrac{1}{2t+1}\cos\alpha\sin^{2t}\alpha \\[1.3mm]
        &- \tfrac{1}{N-2} \cos\alpha\sin^{N-3}\alpha .  
        \end{aligned}  
    \end{aligned}
\end{equation}

\noindent Combining (\ref{odd w}) and (\ref{odd int}), when $N\geq 5$ we have
\begin{equation}
    \label{odd p}
    \begin{aligned}
        p_{N,\alpha} \ &= \ \tfrac{w_{N-2}}{w_{N-1}}\ I_{N-2,\alpha}  \\[1.3mm]
        & \; 
        \begin{aligned}
            = \ \tfrac{(N-2)!!}{2(N-3)!!} \big [&\tfrac{(N-3)!!}{(N-2)!!}\,(1-\cos\alpha )- \tfrac{1}{N-2} \cos\alpha\sin^{N-3}\alpha \\
        &-\sum\limits_{t=1}^{(N-5)/2}\tfrac{(N-3)(N-5)\cdots(2t+2)}{(N-2)(N- 4)\cdots(2t+3)}\tfrac{1}{2t+1}\cos\alpha\sin^{2t}\alpha  \big]
        \end{aligned} \\[1.3mm]
        & = \ \frac{1}{2} \big[ 1-\cos\alpha\big(1+\sum_{t=1}^{(N-3)/2} \tfrac{(2t-1)!!}{(2t)!!}\sin^{2t}\alpha\big)\big]
    \end{aligned}
\end{equation}
\noindent and the proof of the first conclusion is completed. \\
\indent To prove the second conclusion, we firstly consider the condition when $N$ is odd. 
Considering the Taylor series of $h(y)=\frac{1}{\sqrt{1-y}}$ at the point $y_0=0$ and the expansion order $N_1=\frac{N-3}{2}$ with the Lagrange remainder $R_{N_1}(y)$, we have
\begin{equation*}
\label{taylor at 0}
    \begin{aligned}
        h(y)\ &=\ T_{N_1}(h,y_0;y)+R_{N_1}(y) \\
        &= \ \sum_{t=0}^{N_1} \frac{h^{(t)}(y_0)}{t!}(y-y_0)^t+\frac{h^{(N_1+1)}(\xi)}{(N_1+1)!}(y-y_0)^{N_1+1} \quad \big(\xi\in[y_0,y]   \big ) \\
        & =\ 1+\sum_{t=1}^{N_1} \frac{(2t-1)!!}{(2t)!!} y^t+\frac{(2N_1+1)!!}{(2N_1+2)!!}(1-\xi)^{-\frac{3}{2}-N_1} y^{N_1+1} \quad \big( \xi\in[0,y]  \big )  \\
        & =\ 1+\sum_{t=1}^{(N-3)/2} \frac{(2t-1)!!}{(2t)!!} y^t +\frac{(N-2)!!}{(N-1)!!}(1-\xi)^{-\frac{N}{2}} y^{(N-1)/2} \quad \big(\xi\in[0,y]  \big) .
    \end{aligned}
\end{equation*}
Substituting $y$ with $\sin^2\alpha$ and we get
{\small \begin{equation*}
\label{taylor substitute}
\begin{array}{ll}
    \displaystyle \frac{1}{\cos\alpha} = h(\sin^2\alpha)  = 1+\sum_{t=1}^{(N-3)/2} \tfrac{(2t-1)!!}{(2t)!!} \sin^{2t}\alpha  +\tfrac{(N-2)!!}{(N-1)!!}(1-\xi)^{-\frac{N}{2}} \sin^{N-1}\alpha  .\\
\end{array}
\end{equation*}}
for some $\xi\in[0,\sin^2\alpha]$. Thus 
\begin{equation}
\label{ineq 1}
\begin{aligned}
& 1+\sum_{t=1}^{(N-3)/2}\frac{(2t-1)!!}{(2t)!!}\, \sin^{2t}\alpha \\
    =\ & \frac{1}{\cos\alpha}-\frac{(N-2)!!}{(N-1)!!}(1-\xi)^{-\frac{N}{2}} \sin^{N-1}\alpha \quad \big(\xi\in[0,\sin^2\alpha] \,\big)  \\[1.3mm]
    \leq\ &   \frac{1}{\cos\alpha}-\frac{(N-2)!!}{(N-1)!!} \,\sin^{N-1}\alpha  .
\end{aligned}
\end{equation} 

\noindent Substituting (\ref{ineq 1}) into (\ref{odd p}), we have
\begin{equation}
\label{ineq p odd}
\begin{aligned}
    p_{N,\alpha}\ = & \  \frac{1}{2} \,\Big[ 1-\cos\alpha\Big(1+\sum_{t=1}^{(N-3)/2} \frac{(2t-1)!!}{(2t)!!}\sin^{2t}\alpha\Big)\Big]  \\  
    \geq & \ \frac{1}{2}\, \Big [ 1-\cos\alpha \Big (\frac{1}{\cos\alpha}-\frac{(N-2)!!}{(N-1)!!} \, \sin^{N-1}\alpha \Big ) \Big ]  \\[1.3mm]
    = & \  \frac{(N-2)!!}{2(N-1)!!}\,  \cos\alpha\sin^{N-1}\alpha.
\end{aligned}
\end{equation}

\noindent Thus the proof of the second conclusion is completed for any odd $N$. When $N$ is even, the conclusion holds since 
\begin{equation}
\label{ineq p even}
    \begin{aligned}
        p_{N,\alpha} =\ & \frac{ w_{N-2}}{w_{N-1}}\int_0^{\alpha} \sin^{N-2}\theta \, d\theta  \\
        \geq \ &  \frac{w_{N-2}}{w_{N-1}}\int_0^{\alpha} \sin^{N-1}\theta \, d\theta  \\[1.3mm]
        = \ & \frac{w_{N-2}/w_{N-1}}{w_{N-1}/w_{N}\quad}\,  p_{N+1,\alpha} \\[1.3mm]
        \geq \ & \frac{(N-2)!!/\pi(N-3)!!}{(N-1)!!/2(N-2)!!} \, \frac{(N-1)!!}{2\cdot N!!}\, \cos\alpha\sin^{N}\alpha   \\[1.3mm]
        = \ & \frac{(N-2)!!}{\pi N\, (N-3)!!} \, \cos\alpha\sin^{N}\alpha  . \\
    \end{aligned}
\end{equation}
The second inequality above holds by using (\ref{even w}), (\ref{odd w}) and (\ref{ineq p odd}). Therefore the proof of the second conclusion is completed.  
\end{proof}

\begin{figure}
    \centering
    \includegraphics[width=.9\textwidth]{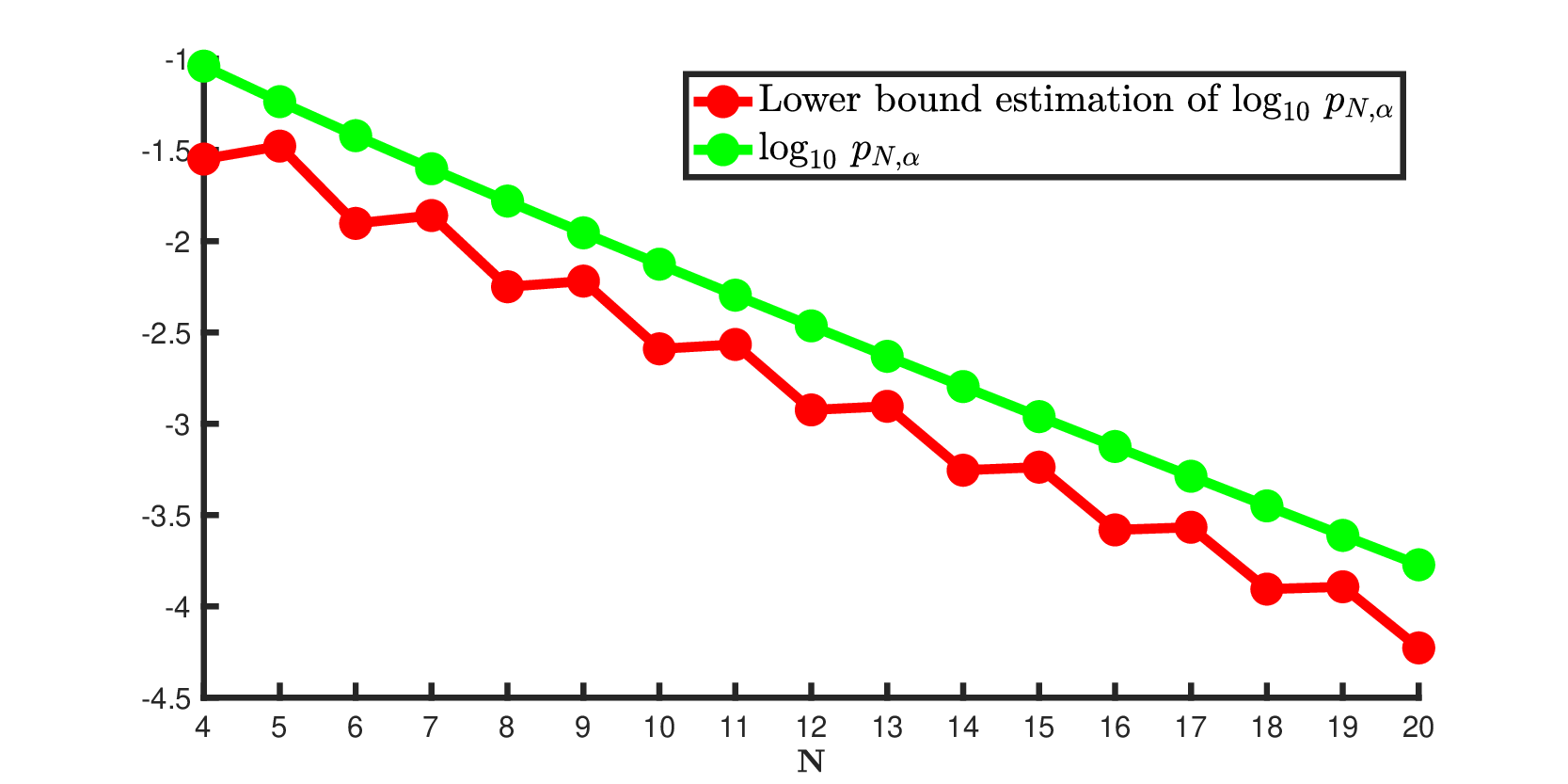}
    \caption{The explicit value and the lower bound of $\log_{10} p_{N,\alpha}$ given in Theorem \ref{theo p} ($\alpha=\frac{\pi}{4}$)}
    \label{fig appro}
\end{figure}

\begin{remark}
Fig. \ref{fig appro} shows the trend of the explicit value and the lower bound of $\log_{10} p_{N,\alpha}$ as $N$ grows ($\alpha=\frac{\pi}{4}$ fixed), using the expressions stated in Theorem \ref{theo p}. The results are shown by their common logarithms. Both curves exhibit linear decays as $N$ grows, illustrating an exponential decay of $p_{N,\alpha}$ on $N$. The lower bound curve is a polygonal line due to its dependence on the parity of $N$. When $N$ is odd, the lower bound is derived directly and is closer to the explicit value; when $N$ is even, the lower bound is deducted from the result of $N+1$, and is therefore farther from the explicit value. 
\end{remark}

 In the following theorem, we establish the lower bound on the probability that DFDS outputs an $\epsilon-$optimal solution, which depends on $p_{N,\alpha}$. This result directly guarantees the convergence in probability of DFDS.

\begin{theorem}
\label{theo p_success}
   DFDS outputs an $\epsilon-$optimal point in $\mathcal{X}$ with probability at least $1-(1-p_{N,\alpha})^M$.
\end{theorem}

\begin{proof}
 It suffices to consider the case when $x_0\notin \mathcal{X}^*_{R_\epsilon}$. In this case, define the \textit{success event} $S$ as the event that the directional random search phase succeeds in finding an $x_k \in \mathcal{X}^*_{R_\epsilon}$ for some $k\geq0$. Following Lemma \ref{lemma optimality}, the sequence of all the iteration points generated in the directional random search phase, denoted by $\{x_k\}_{k=0}^K$, is a sequence of random points satisfying $ K\leq \big\lfloor\frac{f(x_0)-f^*}{\epsilon/3}\big\rfloor$.

Define the random variable $X$ as the last point in $\{x_k\}_{k=0}^K$ that belongs to $\mathcal{X}_{R_\epsilon} \backslash \mathcal{X}^*_{R_\epsilon}$. Then either $X=X_K$ is the terminal point of the directional random search phase, or $X$ is not the last point of the sequence and thus $\{x_k\}_{k=0}^K \cap\mathcal{X}^*_{R_\epsilon}\neq \emptyset$, i.e., $S$ occurs. Following Theorem 1, \vspace{-1ex}
\begin{equation*}
    \begin{aligned}
       P(S)
        = &\; E_{\mathcal{X}_{R_\epsilon} \backslash \mathcal{X}^*_{R_\epsilon}}[P(S|X)] \\
        = &\; 1 - E_{\mathcal{X}_{R_\epsilon} \backslash \mathcal{X}^*_{R_\epsilon}}[P(\text{the directional random search phase terminates at }X|X)] \\ 
        = &\;  1-E[(1-p_{X})^M]\\ 
        \geq &\;   1-(1-p_{N,\alpha})^M .
    \end{aligned} \vspace{-1ex}
\end{equation*}
By Lemma \ref{lemma optimality}, DFDS outputs an $\epsilon-$optimal point in $\mathcal{X}$ with probability at least $P(S)$, and thus the proof is completed.
\end{proof}

By Theorem \ref{theo p}, $ p_{N,\alpha}>0$ when $N$ is fixed. As $M$ (the number of sampling directions) tends to infinity, the probability of outputting an $\epsilon-$optimal point in $\mathcal{X}$ approaches $1$. In other words, DFDS converges in probability to an $\epsilon-$optimal point. 

Having established the convergence in probability, we now turn to analyzing the computational complexity of DFDS. This is measured by the expected number of sampling directions/function evaluations required to attain an $\epsilon-$accuracy. The main result is provided in the following theorem.

\begin{theorem}
    The computational complexity of DFDS is characterized by the following statements. 
    
    1. Define $M_{\epsilon}$ as the number of directions generated in total when an $\frac{2\epsilon}{3}-$optimal iteration point $x_k$ (i.e., $x_k\in\mathcal{X}^*_{R_\epsilon}$) is obtained for the first time. Then the expected value of $M_{\epsilon}$ (denoted as $E[M_{\epsilon}]$) satisfies that when $N$ is even,
    \begin{equation*}
    \lim_{N\rightarrow+\infty} \frac{E[M_{\epsilon}]}{ \sqrt{ N} \big/ \sin^{N}\alpha} \leq \frac{\sqrt{2\pi}}{\cos\alpha} \,\Big\lfloor \frac{f(x_0)-f^*}{\epsilon/3}\Big\rfloor;
    \end{equation*} 
     and when $N$ is odd,
    \begin{equation*}
        \lim_{N\rightarrow+\infty} \frac{E[M_{\epsilon}]}{ \sqrt{ N} \big/ \sin^{N-1}\alpha} \leq  \frac{\sqrt{2\pi}}{\cos\alpha}\, \Big\lfloor \frac{f(x_0)-f^*}{\epsilon/3}\Big\rfloor.
       \end{equation*}
     Furthermore, it is concluded that $E[M_{\epsilon}]$ is at most $O( (\frac{2(D_0+R_\epsilon)}{\sqrt{3}R_\epsilon}) ^{N}\,\sqrt{N}\,$\scalebox{1.25}{$\frac{1}{\epsilon}$}$)$. 
    
    2. To attain an $\epsilon$-optimal solution in $\mathcal{X}$, the expected number of function evaluations in DFDS is at most $O((\frac{2(D_0+R_\epsilon)}{\sqrt{3}R_\epsilon})^{N}\,\frac{D_0}{R_\epsilon}\,\sqrt{N}\,$\scalebox{1.25}{$\frac{1}{\epsilon}$}$)$.
\end{theorem}

\begin{proof}
    It suffices to consider the case when $x_0\notin \mathcal{X}^*_{R_\epsilon}$. Let $K_{hit}$ be the \textit{first hitting time} the iteration point enters $\mathcal{X}^*_{R_\epsilon}$, i.e., $K_{hit}=\inf\{k>0:x_k\in\mathcal{X}^*_{R_\epsilon}\}$. It is clear that \vspace{-1ex}
    \begin{equation*}
    \label{ineq K bound}
       K_{hit}\leq\Big\lfloor \frac{f(x_0)-f^*}{\epsilon/3}\Big\rfloor. \vspace{-1ex} 
    \end{equation*}
    
    At each $x_k$ for $0\leq k\leq K_{hit}-1$ ($x_k$ is a random variable), let $M_\epsilon^k$ be the number of directions sampled to obtain $x_{k+1}$ (an $\frac{\epsilon}{3}-$better point than $x_k$). Conditional on $x_k\notin \mathcal{X}^*_{R_\epsilon}$, $M_\epsilon^k$ follows a geometric distribution with success probability denoted by $p_{x_k}$ and $p_{x_k}\geq p_{N,\alpha}$ based on former analysis. Thus for any $0\leq k\leq K_{hit}-1$, \vspace{-1ex}
    \begin{equation*}
    \label{ineq E(Mk)}
        E[M_\epsilon^k] = E[E[M_\epsilon^k\,|\,x_k]]=E[\,1\, / \, p_{x_k}] \leq 1\, / \, p_{N,\alpha}. \vspace{-1ex}
    \end{equation*}
    Therefore, the expected number of directions generated in total is bounded by \vspace{-1ex}
\begin{equation*} 
    \begin{aligned}
        E[M_\epsilon] = E\big[\sum\limits_{k=0}^{K_{hit}-1} M_\epsilon^k \big] \leq \Big\lfloor \frac{f(x_0)-f^*}{\epsilon/3}\Big\rfloor \, \frac{1}{p_{N,\alpha}} .
    \end{aligned} \vspace{-1ex}
\end{equation*}
    Following the expression of $p_{N,\alpha}$ in Theorem \ref{theo p}, the discussions are separated by the parity of $N$.

Condition 1: $N$ is even.
    
Under this condition, following (\ref{ineq p even}) we have 
\begin{equation*}
\small
\label{even M}
\begin{aligned}
 \lim_{N\rightarrow+\infty} \ \frac{\sin^N\alpha}{\sqrt{ N}} \, \frac{E[M_\epsilon]}{\Big\lfloor \frac{f(x_0)-f^*}{\epsilon/3}\Big\rfloor}  & \leq \ \lim_{N\rightarrow+\infty} \ \frac{\sin^N\alpha}{\sqrt{ N}} \, \frac{1}{p_{N,\alpha}}  \\
    & \leq \  \lim_{N\rightarrow+\infty} \  \frac{1}{\sqrt{ N}}\,\frac{1}{\cos\alpha} \, \frac{\pi N\, (N-3)!!}{(N-2)!!}   \\
    & = \ \frac{\sqrt{2\pi}}{\cos\alpha}\, \lim_{N\rightarrow+\infty} \ \frac{\sqrt{\pi N}(N-3)!!}{\sqrt{2}(N-2)!!}  
    \\
    & = \ \frac{\sqrt{2\pi}}{\cos\alpha} .
\end{aligned} \vspace{0.5ex}
\end{equation*}
For the establishment of the last equality, the limits \vspace{-1ex}
\begin{equation*}
    \lim_{N\rightarrow+\infty} \ \frac{\sqrt{\pi N}(N-3)!!}{\sqrt{2}(N-2)!!}=1  \vspace{-1ex}
\end{equation*}
is deducted from the Wallis formula \cite{Wastlund2007} \vspace{-1ex}
\begin{equation}
\label{eq wallis}
    \displaystyle\lim\limits_{j\rightarrow+\infty}\frac{1}{2j+1}\Big[\frac{(2j)!!}{(2j-1)!!}\Big]^2 = \frac{\pi}{2}. \vspace{-1ex}
\end{equation} 

Condition 2: $N$ is odd.
 
Under this condition, following (\ref{ineq p odd}) we have
\begin{equation*}
\small
\label{odd M}
\begin{aligned}
\lim_{N\rightarrow+\infty} \frac{\sin^{N-1}\alpha}{\sqrt{ N}} \, \frac{E[M_\epsilon]}{\Big\lfloor \frac{f(x_0)-f^*}{\epsilon/3}\Big\rfloor} & \leq \ \lim_{N\rightarrow+\infty} \frac{\sin^{N-1}\alpha}{\sqrt{ N}} \, \frac{1}{p_{N,\alpha}}  \\
    & \leq \ \lim_{N\rightarrow+\infty} \frac{1}{\sqrt{N}}\, \frac{1}{\cos\alpha } \, \frac{2(N-1)!!}{(N-2)!!}\,   \\
    & 
    = \ \frac{\sqrt{2\pi}}{\cos\alpha}  \, \lim_{N\rightarrow+\infty} \frac{\sqrt{2}(N-1)!!}{\sqrt{\pi N}(N-2)!!}         \\
    & = \ \frac{\sqrt{2\pi}}{\cos\alpha}  .
\end{aligned} \vspace{0.5ex}
\end{equation*}
For the establishment of the last equality, the limits \vspace{-1ex}
\begin{equation*}
    \lim_{N\rightarrow+\infty} \frac{\sqrt{2}(N-1)!!}{\sqrt{\pi N}(N-2)!!}  = 1 \vspace{-1ex}
\end{equation*}
is deducted from the Wallis formula (\ref{eq wallis}). \vspace{1ex}

By Lemma \ref{lemm alpha and D}, $\frac{1}{\sin\alpha}\leq \frac{2(D_0+R_\epsilon)}{\sqrt{3}R_\epsilon}$, and $\frac{1}{\cos\alpha}<\frac{1}{\cos(\pi/3)}=2$ is a constant. Thus it is concluded that $E[M_{\epsilon}]$ is at most $O\big( (\frac{2(D_0+R_\epsilon)}{\sqrt{3}R_\epsilon}) ^{N}\,\sqrt{N}\,$\scalebox{1.25}{$\frac{1}{\epsilon}$}$)$ for all $N$ and the first statement is proved.

Along each direction, points are searched at a distance of $R_\epsilon$ inside $\mathcal{X}_{R_\epsilon}$, generating at most $\lfloor\frac{D(\mathcal{X}_{R_\epsilon})} {R_\epsilon}\rfloor\leq\lfloor\frac{D_0}{R_\epsilon}\rfloor+2$ function evaluations. Together with Lemma \ref{lemma optimality}, it is concluded that DFDS finds an $\epsilon-$optimal solution in $\mathcal{X}$ with at most an expected number of $O((\frac{2(D_0+R_\epsilon)}{\sqrt{3}R_\epsilon})^{N}\,\frac{D_0}{R_\epsilon}\,\sqrt{N}\,$\scalebox{1.25}{$\frac{1}{\epsilon}$}$)$ function evaluations.
\end{proof}

Having established results on the convergence and complexity of DFDS, we conclude this section with discussions on $R_\epsilon$. As stated in Section \ref{section preliminaries}, the existence of $R_\epsilon$ is guaranteed for (\ref{problem}). Moreover, an explicit expression for $R_\epsilon$ can be derived when $f$ satisfies additional regularity conditions. \vspace{1ex}

\textbf{1. Lipschitz continuity}

If $f$ is Lipschitz continuous with constant $L$ on some extended set $\mathcal{X}_R$, i.e., \vspace{-1ex}
\begin{equation*}
    |f(x)-f(y)|\leq L\|x-y\|\quad \text{for all } x,y\in\mathcal{X}_R, \vspace{-1ex}
\end{equation*}
then choosing $R_\epsilon=\min\{R,\frac{\epsilon}{3L}\}$ satisfies Assumption \ref{asp R0}.  \vspace{1ex}

For example, if $f(x)=c^Tx$ is linear, then $L=\|c\|$ and $R_\epsilon=\frac{\epsilon}{3\|c\|}$ can be taken.  \vspace{1.2ex}

\textbf{2. Continuous differentiability}

If $f$ is continuously differentiable on some bounded extended set $\mathcal{X}_R$, then $\nabla f$ is bounded on $\mathcal{X}_R$. Let $M>0$ be such that $\|\nabla f(x)\|\leq M$ for all $x\in\mathcal{X}_R$. By the mean value theorem, 
\begin{equation*}
    |f(x)-f(y)|\leq M\|x-y\|\quad \text{for all } x,y\in\mathcal{X}_R, \vspace{1.2ex}
\end{equation*}
and $R_\epsilon = \min\{R,\frac{\epsilon}{3M}\}$ can be taken.  \vspace{1.2ex}

For example, if $f(x)=\frac{1}{2}x^TQx+q^Tx$ is quadratic, then $\nabla f(x)=Q^Tx+q$ is continuous, and the above applies.

\section{Numerical Experiments}
\label{section numerical}

This section compares DFDS against two baseline random search methods: pure random search (PRS) and improving hit-and-run (IHR). Numerical results demonstrate DFDS's dominance in efficiency and solution accuracy, achieving higher success rates in locating global minima under reduced function evaluation budgets across scalable dimensional problems.

\subsection{Background}

The algorithms are implemented in Matlab R2024a on a MacBook Air laptop with an Apple M2 chip and 16GB memory. The error threshold is set as $\epsilon= 10^{-4}$, and in experiments, every solution within this threshold is regarded as global optimal. 

Three random search methods are performed in total, including DFDS, PRS in \cite{Zabinsky2003_PRS} and IHR in \cite{Zabinsky1993}. The baseline selection is motivated by two considerations: PRS represents an implementable sampling method sharing the same exponential complexity class as DFDS; IHR provides a contrasting directional search method, generating the searching sequence by broad coverage exploration, whereas DFDS adopts a depth-first strategy. The parameters included in the algorithms are defined as follows:
\begin{itemize}
  \item \textbf{DFDS}: 
       $R_0$ represents the uniform line search distance; $M_{\text{DFDS}}$ represents the maximum random directions per iteration.
  \item \textbf{IHR}: $M_{\text{IHR}}$ represents the maximum random directions per iteration.
  \item \textbf{PRS}: $n_{\text{PRS}}$ represents the number of sampling points in total.
\end{itemize}

For each random search algorithm, there exists a sampling parameter: $M_{\text{DFDS}}$ in DFDS, $M_{\text{IHR}}$ in IHR and $n_{\text{PRS}}$ in PRS. In our experiments, we implemented a \textit{uniform function evaluation budget} ($n_\mathrm{feval}$) across all algorithms: solution quality metrics are compared after each algorithm executes $n_\mathrm{feval}$ function evaluations. As the computational cost is dominated by function evaluations, this design ensures nearly identical runtimes for all three methods. Three budget levels of $n_\mathrm{feval}$ were specified for each problem, defined as exponential functions of $N$ for variable-dimensional test problems. 

Regarding the parameter $R_0$, we forwent deriving an exact value that satisfies the theoretical assumptions for two reasons: first, computing it is difficult for certain global optimization benchmarks; second, the theoretical value provides a worst-case lower bound that is often overly conservative. Given that $R_0$ conceptually represents a radius, we selected it empirically to be proportional to $\sqrt{N}$. Such choice is justified by its consistent and robust performance across problem instances.

Exact forms of parameters are given in Table \ref{table parameter}.

\begin{table}[ht]
\onehalfspacing
\captionsetup{justification=raggedright, singlelinecheck=false}
\caption{Parameter settings for benchmark problems} \vspace{-1ex}
\label{table parameter}
\begin{tabular}{@{}l c c c c c c c@{}}
\toprule
\headerVcenter[5pt]{2}{\textbf{Problem}} & 
\multicolumn{3}{c}{\bm{$n_\mathrm{feval}$}} & 
\headerVcenter[5pt]{2}{\textbf{DFDS} \bm{$R_0$}} \\
\cmidrule(lr){2-4} 
& \textbf{Low} & \textbf{Medium} & \textbf{High} &  \\
\midrule
Six-Hump Camel & $125$ & $250$ & $500$ & $0.5$  \\
Goldstein Price & $125$ & $250$ & $500$ & $0.2$\\
Ackley / Levy & $125\times2^N$ & $250\times2^N$ & $500\times2^N$ & $\frac{1}{2\sqrt{2}}\sqrt{N}$  \\
Alpine & $625\times2^N$ & $1250\times2^N$ & $2500\times2^N$ & $\frac{1}{2\sqrt{2}}\sqrt{N}$  \\
\bottomrule
\end{tabular}

\small{
\begin{tabular}{@{}p{\textwidth}@{}}
\bm{$n_\mathrm{feval}$} denotes the function evaluation budget per run (uniform across algorithms). \\
\bm{$R_0$} denotes the line search distance (DFDS only). 
\end{tabular}
}
\end{table}

The algorithms were evaluated on well-established baseline test problems drawn from \cite{Jamil2013}. Variable-dimensional functions were specifically selected to assess performance evolution across dimensions ($N$).

1. Six-Hump Camel Function
\begin{equation*}
    f(x)= 4x_1^2-2.1x_1^4+\frac{1}{3}x_1^6+x_1x_2-4x_2^2+4x_2^4
\end{equation*}
subject to $-5\leq x_1,x_2 \leq 5$.
The global minimum $f^* = -1.0316$ is attained at $x_1^*=\begin{pmatrix}
    -0.0898 & 0.7126\end{pmatrix}^T$ and $x_2^*=\begin{pmatrix} 0.0898 & -0.7126 \end{pmatrix}^T$.\vspace{1.5ex}    
    
2. Goldstein Price Function
\begin{equation*}
\begin{aligned}
    f(x)= &[1+(x_1+x_2+1)^2(19-14x_1+3x_1^2-14x_2+6x_1x_2+3x_2^2)] \\
    &\times [30+(2x_1-3x_2)^2(18-32x_1+12x_1^2+48x_2-36x_1x_2+27x_2^2)]  
\end{aligned}
\end{equation*}
subject to $-2\leq x_1,x_2 \leq 2$. The global minimum $f^* = 3$ is attained at $x^*=\begin{pmatrix} 0 & -1\end{pmatrix}^T$.  \vspace{1.5ex}

3. Ackley Function:
\begin{equation*}
    f(x)=-20\cdot e^{-0.2\sqrt{\frac{1}{N}\sum_{i=1}^Nx_i^2} \,} - e^{\frac{1}{N}\sum_{i=1}^N\cos(2\pi x_i)}+20+\exp(1)
\end{equation*}
subject to $-10\leq x_i \leq 10, \; i=1,...,N.$ The global minimum $f^* = 0$ is attained at $x^*=\begin{pmatrix} 0 & 0 & \cdots & 0
\end{pmatrix}^T$.  \vspace{1.5ex}
    
4. Levy Function:
\begin{equation*}
    f(x) = \sin^2(\pi y_1)+\sum\limits_{i=1}^{N-1}(y_i-1)^2[1+10\sin^2(\pi y_i+1)]+(y_N-1)^2[1+\sin^2(2\pi y_N)] 
\end{equation*}
subject to $y_i=1+\frac{x_i-1}{4}$, $-10\leq x_i \leq 10, \; i=1,...,N.$ The global minimum $f^* = 0$ is attained at $x^*=\begin{pmatrix} 1 & 1 & \cdots &1
\end{pmatrix}^T$.  \vspace{1.5ex}

5. Alpine Function:
\begin{equation*}
    f(x) = -\prod\limits_{i=1}^N \sqrt{x_i}\sin(x_i)
\end{equation*}
subject to $0\leq x_i \leq 10, \; i=1,...,N.$ The global minimum $f^* \approx -2.808 ^ N$ is attained at $x^*\approx\begin{pmatrix} 7.917 & 7.917 & \cdots & 7.917 \end{pmatrix}^T$.  \vspace{1.5ex}

We evaluated the algorithms on low- and moderate-dimensional instances of the above five benchmark functions. The first two functions are of fixed dimensionality, while the latter three are scalable to arbitrary dimensions. For each instance, we conducted 10 independent runs, where the initial points were uniformly sampled from $\mathcal{X}$ (note that PRS does not require an initial point). To refine the solution precision, a local search step was applied to the final point obtained by each algorithm. 

The performance was assessed based on two solution quality metrics: the success rate ($\text{SR}$) of locating $\epsilon-$optimal solutions ($\epsilon=10^{-4}$), and the minimal function value found ($f_{\mathrm{best}}$). Given that runtimes are controlled to be nearly identical, we omitted the runtime indicator and focused solely on solution quality. 

 For variable-dimensional problems, a stopping criterion was applied: the increase in dimensionality was halted once the accuracy of at least two algorithms falls below $50\%$.

\begin{table}[ht]
\captionsetup{justification=raggedright,singlelinecheck=false}
\caption{Algorithm performance on functions with fixed dimensionality} 
\small
\vspace{-0.5em}
\label{table N=2}
\begin{tabular}{@{}c c c c c c c c c c c@{}}
\toprule
\textbf{Problem} & \textbf{Budget} & \large\bm{$n_\mathrm{feval}$} & \textbf{DFDS SR} & \textbf{IHR SR} & \textbf{PRS SR}\\
\midrule

\vcenteredmultirow{3}{\shortstack{Goldstein Price}} & \textbf{Low} & 125 & 70\% & 100\% & 90\% \\
\cmidrule(lr){2-6}
   & \textbf{Medium} & 250 & 80\% & 100\% & 100\% \\
\cmidrule(lr){2-6}
   & \textbf{High} & 500 & 100\% & 100\% & 100\% \\
\midrule

\vcenteredmultirow{3}{Six-Hump Camel} & \textbf{Low} & 125 & 40\% & 100\% & 100\% \\
\cmidrule(lr){2-6}
  & \textbf{Medium} & 250 & 100\% & 100\% & 100\% \\
\cmidrule(lr){2-6}
  & \textbf{High} & 500 & 100\% & 100\% & 100\% \\
\midrule

\vcenteredmultirow{3}{Ackley ($N=2$)} & \textbf{Low} & 500 & 70\% & 100\% & 80\% \\
\cmidrule(lr){2-6}
    & \textbf{Medium} & 1,000 & 100\% & 100\% & 80\% \\
\cmidrule(lr){2-6}
    & \textbf{High} & 2,000 & 100\% & 100\% & 100\% \\
\midrule

\vcenteredmultirow{3}{Levy ($N=2$)} & \textbf{Low} & 500 & 80\% & 100\% & 80\% \\
\cmidrule(lr){2-6}
  & \textbf{Medium} & 1,000 & 100\% & 100\% & 80\% \\
\cmidrule(lr){2-6}
  & \textbf{High} & 2,000 & 100\% & 100\% & 100\% \\
  
\bottomrule
\end{tabular}

\small{
\begin{tabular}{@{}p{\textwidth}@{}}
\textbf{Budget} denotes the function evaluation budget level (Low/Medium/High), 10 independent runs are executed per problem at each budget level. \\
\bm{$n_\mathrm{feval}$} denotes the function evaluation budget per run. \\
\textbf{SR} denotes success rate (\%) of locating $\epsilon-$optimal solutions ($\epsilon = 10^{-4}$). 
\end{tabular}}
\end{table}

\begin{table}[ht]
\captionsetup{justification=raggedright,singlelinecheck=false}
\caption{Algorithm performance on the Ackley Function} 
\vspace{-0.5em}
\label{table Ackley}
\small
\begin{tabular}{@{}c c c c c c c c c c c@{}}
\toprule
\headerVcenter[4pt]{2}{\bm{$N$}} & \headerVcenter[4pt]{2}{\textbf{Budget}} & \headerVcenter[4pt]{2}{\large\bm{$n_\mathrm{feval}$}} & \multicolumn{2}{c}{\textbf{DFDS}} & \multicolumn{2}{c}{\textbf{IHR}} & \multicolumn{2}{c}{\textbf{PRS}} \\
\cmidrule(lr){4-5} \cmidrule(lr){6-7} \cmidrule(lr){8-9}
 & & & \textbf{SR} & \normalsize\bm{$f_\mathrm{best}$} & \textbf{SR} & \normalsize\bm{$f_\mathrm{best}$} & \textbf{SR} & \normalsize\bm{$f_\mathrm{best}$} \\
\midrule

\vcenteredmultirow{3}{5} & \textbf{Low} & 4,000 & 100\% & $3.55$e-15 & 60\% & $3.55 $e-15 & 20\% & $3.55$e-15 \\
\cmidrule(lr){2-9}
       & \textbf{Medium} & 8,000 & 100\% & $3.55$e-15 & 60\% & $3.55 $e-15 & 30\% & $3.55 $e-15 \\
\cmidrule(lr){2-9}
       & \textbf{High} & 16,000 & 100\% & $3.55 $e-15 & 90\% & $3.55 $e-15 & 40\% & $3.55 $e-15 \\
\midrule

\vcenteredmultirow{3}{6} & \textbf{Low} & 8,000 & 80\% & $3.55 $e-15 & 60\% & $3.55 $e-15 & 0\% & 1.50 \\
\cmidrule(lr){2-9}
       & \textbf{Medium} & 16,000 & 100\% & $3.55 $e-15 & 70\% & $3.55 $e-15 & 10\% & $7.11 $e-15 \\
\cmidrule(lr){2-9}
       & \textbf{High} & 32,000 & 100\% & $3.55 $e-15 & 80\% & $3.55 $e-15 & 10\% & $7.11 $e-15  \\
\midrule

\vcenteredmultirow{3}{7} & \textbf{Low} & 16,000 & 90\% & $3.55 $e-15 & 90\% & $3.55 $e-15 & 10\% & $7.11 $e-15 \\
\cmidrule(lr){2-9}
       & \textbf{Medium} & 32,000 & 100\% & $3.55 $e-15 & 100\% & $3.55 $e-15 & 10\% & $7.11 $e-15 \\
\cmidrule(lr){2-9}
       & \textbf{High} & 64,000 & 100\% & $3.55 $e-15 & 100\% & $3.55 $e-15 & 10\% & $7.11 $e-15 \\
\midrule

\vcenteredmultirow{3}{8} & \textbf{Low} & 32,000 & 100\% & $7.11 $e-15 & 60\% & $1.50 $e-6 & 0\% & 2.25 \\
\cmidrule(lr){2-9}
       & \textbf{Medium} & 64,000 & 100\% & $3.55 $e-15 & 70\% & $3.72 $e-6 & 0\% & 2.25 \\
\cmidrule(lr){2-9}
       & \textbf{High} & 128,000 & 100\% & $3.55 $e-15 & 80\% & $1.40 $e-6 & 0\% & 2.25 \\
\midrule

\vcenteredmultirow{3}{9} & \textbf{Low} & 64,000 & 80\% & $3.55 $e-15 & 60\% & $1.61 $e-6 & 10\% & $7.11 $e-15 \\
\cmidrule(lr){2-9}
       & \textbf{Medium} & 128,000 & 90\% & $3.55 $e-15 & 70\% & $1.04 $e-6 & 10\% & $7.11 $e-15 \\
\cmidrule(lr){2-9}
       & \textbf{High} & 256,000 & 90\% & $3.55 $e-15 & 90\% & $2.43 $e-7 & 20\% & $7.11 $e-15 \\
\midrule

\vcenteredmultirow{3}{10} & \textbf{Low} & 128,000 & 70\% & $7.11 $e-15 & 60\% & $7.98 $e-8 & 0\% & 1.65 \\
\cmidrule(lr){2-9}
        & \textbf{Medium} & 256,000 & 90\% & $3.55 $e-15 & 70\% & $2.13 $e-7 & 0\% & 1.65 \\
\cmidrule(lr){2-9}
        & \textbf{High} & 512,000 & 90\% & $3.55 $e-15 & 80\% & $2.28 $e-6 & 0\% & 1.65 \\
\midrule

\vcenteredmultirow{3}{11} & \textbf{Low} & 256,000 & 90\% & $7.11 $e-15 & 20\% & $2.01 $e-5 & 0\% & 1.10 \\
\cmidrule(lr){2-9}
        & \textbf{Medium} & 512,000 & 100\% & $7.11 $e-15 & 20\% & $7.15 $e-6 & 10\% & $7.11 $e-15 \\
\cmidrule(lr){2-9}
        & \textbf{High} & 1,024,000 & 100\% & $7.11 $e-15 & 60\% & $6.84 $e-7 & 10\% & $7.11 $e-15 \\
\midrule

\vcenteredmultirow{3}{12} & \textbf{Low} & 512,000 & 80\% & $7.11 $e-15 & 30\% & $2.42 $e-6 & 0\% & 2.12 \\
\cmidrule(lr){2-9}
        & \textbf{Medium} & 1,024,000 & 90\% & $7.11 $e-15 & 40\% & $7.36 $e-7 & 0\% & 1.84 \\
\cmidrule(lr){2-9}
        & \textbf{High} & 2,048,000 & 90\% & $7.11 $e-15 & 50\% & $5.20 $e-7 & 0\% & 1.84 \\
\bottomrule
\end{tabular}
\footnotesize{
\begin{tabular}{@{}p{\textwidth}@{}}
\textbf{Budget} denotes the function evaluation budget level (Low/Medium/High), 10 independent runs are executed per instance at each budget level. \\
\bm{$n_\mathrm{feval}$} denotes the function evaluation budget per run. \\
\textbf{SR} denotes success rate (\%) of locating $\epsilon-$optimal solutions ($\epsilon = 10^{-4}$). \\
\bm{$f_\mathrm{best}$} denotes the minimal function value found in each instance.
\end{tabular}}
\end{table}

\begin{table}[ht]
\captionsetup{justification=raggedright,singlelinecheck=false}
\caption{Algorithm performance on the Levy Function} 
\vspace{-0.5em}
\label{table Levy}
\small
\begin{tabular}{@{}c c c c c c c c c c@{}}
\toprule
\headerVcenter[4pt]{2}{\bm{$N$}} & \headerVcenter[4pt]{2}{\textbf{Budget}} & \headerVcenter[4pt]{2}{\large\bm{$n_\mathrm{feval}$}} & \multicolumn{2}{c}{\textbf{DFDS}} & \multicolumn{2}{c}{\textbf{IHR}} & \multicolumn{2}{c}{\textbf{PRS}} \\
\cmidrule(lr){4-5} \cmidrule(lr){6-7} \cmidrule(lr){8-9}
 & & & \textbf{SR} & \normalsize\bm{$f_\mathrm{best}$} & \textbf{SR} & \normalsize\bm{$f_\mathrm{best}$} & \textbf{SR} & \normalsize\bm{$f_\mathrm{best}$} \\
 
\midrule
\vcenteredmultirow{3}{5} & \textbf{Low} & 4,000 & 90\% & $1.16 $e-10 & 100\% & $2.74 $e-11 & 20\% & $1.08 $e-10 \\
\cmidrule(lr){2-9}
  & \textbf{Medium} & 8,000 & 100\% & $4.26 $e-11 & 100\% & $3.69 $e-11 & 30\% & $4.99 $e-11 \\
\cmidrule(lr){2-9}
  & \textbf{High} & 16,000 & 100\% & $2.47 $e-11 & 100\% & $1.19 $e-10 & 40\% & $4.99 $e-11 \\
\midrule

\vcenteredmultirow{3}{6} & \textbf{Low} & 8,000 & 100\% & $1.35 $e-11 & 80\% & $2.44 $e-11 & 10\% & $4.81 $e-11 \\
\cmidrule(lr){2-9}
  & \textbf{Medium} & 16,000 & 100\% & $3.79 $e-11 & 80\% & $3.57 $e-11 & 10\% & $4.81 $e-11 \\
\cmidrule(lr){2-9}
  & \textbf{High} & 32,000 & 100\% & $3.79 $e-11 & 100\% & $4.32 $e-11 & 20\% & $5.51 $e-12 \\
\midrule

\vcenteredmultirow{3}{7} & \textbf{Low} & 16,000 & 90\% & $3.52 $e-11 & 80\% & $8.89 $e-11 & 0\% & 0.04 \\
\cmidrule(lr){2-9}
  & \textbf{Medium} & 32,000 & 90\% & $9.50 $e-11 & 100\% & $7.23 $e-11 & 0\% & 0.04 \\
\cmidrule(lr){2-9}
  & \textbf{High} & 64,000 & 100\% & $9.50 $e-11 & 100\% & $5.46 $e-11 & 10\% & $2.91 $e-9 \\
\midrule

\vcenteredmultirow{3}{8} & \textbf{Low} & 32,000 & 80\% & $7.34 $e-11 & 60\% & $2.27 $e-11 & 0\% &  0.04\\
\cmidrule(lr){2-9}
  & \textbf{Medium} & 64,000 & 90\% & $2.81 $e-11 & 100\% & $1.78 $e-11 & 0\% & 0.04 \\
\cmidrule(lr){2-9}
  & \textbf{High} & 128,000 & 100\% & $3.07 $e-11 & 100\% & $2.27 $e-11 & 10\% & $1.49 $e-10 \\
\midrule

\vcenteredmultirow{3}{9} & \textbf{Low} & 64,000 & 90\% & $7.82 $e-12 & 70\% & $7.34 $e-12 & 0\% & 0.10 \\
\cmidrule(lr){2-9}
  & \textbf{Medium} & 128,000 & 100\% & $7.82 $e-12 & 70\% & $3.18 $e-11 & 0\% & 0.10 \\
\cmidrule(lr){2-9}
  & \textbf{High} & 256,000 & 100\% & $2.96 $e-11 & 80\% & $1.34 $e-11 & 0\% & 0.09 \\
\midrule

\vcenteredmultirow{3}{10} & \textbf{Low} & 128,000 & 90\% & $5.86 $e-12 & 50\% & $2.12 $e-11 & 0\% & 0.20 \\
\cmidrule(lr){2-9}
   & \textbf{Medium} & 256,000 & 100\% & $5.86 $e-12 & 50\% & $2.12 $e-11 & 0\% & 0.17 \\
\cmidrule(lr){2-9}
   & \textbf{High} & 512,000 & 100\% & $5.86 $e-12 & 80\% & $2.12 $e-11 & 0\% & 0.17 \\
\midrule

\vcenteredmultirow{3}{11} & \textbf{Low} & 256,000 & 50\% & $8.40 $e-11 & 60\% & $1.23 $e-11 & 0\% & 0.13 \\
\cmidrule(lr){2-9}
   & \textbf{Medium} & 512,000 & 90\% & $6.98 $e-12 & 60\% & $1.23 $e-11 & 0\% & 0.13 \\
\cmidrule(lr){2-9}
   & \textbf{High} & 1,024,000 & 100\% & $6.98 $e-12 & 90\% & $1.23 $e-11 & 0\% & 0.05 \\
\midrule

\vcenteredmultirow{3}{12} & \textbf{Low} & 512,000 & 70\% & $1.39 $e-11 & 50\% & $1.74 $e-10 & 0\% & 0.19 \\
\cmidrule(lr){2-9}
   & \textbf{Medium} & 1,024,000 & 90\% & $1.39 $e-11 & 60\% & $1.76 $e-11 & 0\% & 0.19 \\
\cmidrule(lr){2-9}
   & \textbf{High} & 2,048,000 & 100\% & $1.39 $e-11 & 60\% & $1.76 $e-11 & 0\% & 0.05 \\
\midrule

\vcenteredmultirow{3}{13} & \textbf{Low} & 1,024,000 & 80\% & $3.99 $e-11 & 30\% & $1.03 $e-10 & 0\% & 0.35 \\
\cmidrule(lr){2-9}
   & \textbf{Medium} & 2,048,000 & 80\% & $2.36 $e-11 & 40\% & $1.03 $e-10 & 0\% & 0.26 \\
\cmidrule(lr){2-9}
   & \textbf{High} & 4,096,000 & 90\% & $1.40 $e-11 & 70\% & $7.97 $e-12 & 0\% & 0.13 \\
\midrule

\vcenteredmultirow{3}{14} & \textbf{Low} & 2,048,000 & 90\% & $1.29 $e-11 & 40\% & $1.85 $e-11 & 0\% & 0.45 \\
\cmidrule(lr){2-9}
   & \textbf{Medium} & 4,096,000 & 90\% & $1.29 $e-11 & 50\% & $1.85 $e-11 & 0\% & 0.23 \\
\cmidrule(lr){2-9}
   & \textbf{High} & 8,192,000 & 90\% & $1.29 $e-11 & 50\% & $1.85 $e-11 & 0\% & 0.20 \\
\bottomrule
\end{tabular}
\vspace{-0.5em}
\end{table}

\begin{table}[ht]
\captionsetup{justification=raggedright,singlelinecheck=false}
\caption{Algorithm performance on the Alpine Function} 
\vspace{-0.5em}
\label{table Alpine}
\small
\begin{tabular}{@{}c c c c c c c c c c c@{}}
\toprule
\headerVcenter[4pt]{2}{\bm{$N$}} & \headerVcenter[4pt]{2}{\textbf{Budget}} & \headerVcenter[4pt]{2}{\large\bm{$n_\mathrm{feval}$}} & \multicolumn{2}{c}{\textbf{DFDS}} & \multicolumn{2}{c}{\textbf{IHR}} & \multicolumn{2}{c}{\textbf{PRS}} \\
\cmidrule(lr){4-5} \cmidrule(lr){6-7} \cmidrule(lr){8-9}
 & & & \textbf{SR} & \normalsize\bm{$f_\mathrm{best}$} & \textbf{SR} & \normalsize\bm{$f_\mathrm{best}$} & \textbf{SR} & \normalsize\bm{$f_\mathrm{best}$} \\
\midrule
\vcenteredmultirow{3}{2} & \textbf{Low} & 2,500 & 100\% & -7.89 & 100\% & -7.89 & 100\% & -7.89 \\
\cmidrule(lr){2-9}
   & \textbf{Medium} & 5,000 & 100\% & -7.89 & 100\% & -7.89 & 100\% & -7.89 \\
\cmidrule(lr){2-9}
   & \textbf{High} & 10,000 & 100\% & -7.89 & 100\% & -7.89 & 100\% & -7.89 \\
\midrule

\vcenteredmultirow{3}{3} & \textbf{Low} & 5,000 & 100\% & -22.14 & 100\% & -22.14 & 100\% & -22.14 \\
\cmidrule(lr){2-9}
   & \textbf{Medium} & 10,000 & 100\% & -22.14 & 100\% & -22.14 & 100\% & -22.14 \\
\cmidrule(lr){2-9}
   & \textbf{High} & 20,000 & 100\% & -22.14 & 100\% & -22.14 & 100\% & -22.14 \\
\midrule
\vcenteredmultirow{3}{4} & \textbf{Low} & 10,000 & 100\% & -62.18 & 90\% & -62.18 & 100\% & -62.18 \\
\cmidrule(lr){2-9}
   & \textbf{Medium} & 20,000 & 100\% & -62.18 & 100\% & -62.18 & 100\% & -62.18 \\
\cmidrule(lr){2-9}
   & \textbf{High} & 40,000 & 100\% & -62.18 & 100\% & -62.18 & 100\% & -62.18 \\
\midrule

\vcenteredmultirow{3}{5} & \textbf{Low} & 20,000 & 70\% & -174.62 & 70\% & -174.62 & 80\% & -174.62 \\
\cmidrule(lr){2-9}
   & \textbf{Medium} & 40,000 & 70\% & -174.62 & 90\% & -174.62 & 100\% & -174.62 \\
\cmidrule(lr){2-9}
   & \textbf{High} & 80,000 & 100\% & -174.62 & 100\% & -174.62 & 100\% & -174.62 \\
\midrule

\vcenteredmultirow{3}{6} & \textbf{Low} & 40,000 & 40\% & -490.35 & 30\% & -490.35 & 60\% & -490.35 \\
\cmidrule(lr){2-9}
   & \textbf{Medium} & 80,000 & 80\% & -490.35 & 50\% & -490.35 & 60\% & -490.35 \\
\cmidrule(lr){2-9}
   & \textbf{High} & 160,000 & 100\% & -490.35 & 60\% & -490.35 & 90\% & -490.35 \\
\midrule
\vcenteredmultirow{3}{7} & \textbf{Low} & 80,000 & 40\% & $-1.38 $e3 & 40\% & $-1.38 $e3 & 30\% & $-1.38 $e3 \\
\cmidrule(lr){2-9}
   & \textbf{Medium} & 160,000 & 40\% & $-1.38 $e3 & 40\% & $-1.38 $e3 & 30\% & $-1.38 $e3 \\
\cmidrule(lr){2-9}
   & \textbf{High} & 320,000 & 100\% & $-1.38 $e3 & 40\% & $-1.38 $e3 & 60\% & $-1.38 $e3 \\
\midrule
\vcenteredmultirow{3}{8} & \textbf{Low} & 160,000 & 20\% & $-3.87 $e3 & 0\% & $-2.34 $e3 & 20\% & $-3.87 $e3 \\
\cmidrule(lr){2-9}
   & \textbf{Medium} & 320,000 & 20\% & $-3.87 $e3 & 10\% & $-3.87 $e3 & 20\% & $-3.87 $e3 \\
\cmidrule(lr){2-9}
   & \textbf{High} & 640,000 & 60\% & $-3.87 $e3 & 20\% & $-3.87 $e3 & 20\% & $-3.87 $e3 \\
\bottomrule
\end{tabular}
\vspace{-0.5em}
\end{table}

Table \ref{table N=2} presents numerical results for fixed-dimensional problems. All three algorithms achieve $100\%$ solution accuracy under high evaluation budgets, confirming their global optimization capability. IHR dominates in three algorithms by maintaining $100\%$ accuracy across all budget levels. PRS exhibits sporadic failures on Ackley and Levy functions, likely due to their expansive feasible domains. DFDS exhibits a lower $\mathrm{SR}$ on the Goldstein Price and Six-Hump Camel functions under constrained evaluation budgets. This can be attributed to its depth-first line search strategy, resulting in fewer directions explored when the total budget is limited.

Table \ref{table Ackley}-\ref{table Alpine} display numerical results for variable-dimensional instances. Among three algorithms, PRS performs the worst. Its performance declines sharply with increasing $N$, achieving $\mathrm{SR}\leq20\%$ for $N\geq 8$ and failing to locate global minima in Levy/Ackley functions at moderate dimensions. DFDS surpasses IHR by both solution accuracy and quality in Table \ref{table Ackley} (Ackley). DFDS consistently matches or exceeds IHR's $\mathrm{SR}$ across all dimensions and budget levels. Notably, DFDS's medium-budget accuracy remains superior to IHR's high-budget performance at $N\geq 8$, highlighting the efficiency advantage of DFDS. Moreover, DFDS yields $f_{\mathrm{best}}$ values closer to $f^*$ than IHR by orders of magnitude. In Table \ref{table Levy} (Levy), DFDS and IHR both maintain $\mathrm{SR}\geq80\%$ for $N\leq 8$, and DFDS attains $100\%$ $\text{SR}$ under high budgets at moderate dimensions, while IHR exhibits occasional failures. In Table \ref{table Alpine} (Alpine), DFDS and IHR both maintain $\text{SR}\geq 70\%$ for $N\leq 5$, while DFDS achieves $\geq 30\%$ higher SR for $6\leq N\leq 9$ under high budgets.

 In summary, DFDS demonstrates superior accuracy: maintaining $\mathrm{SR}\geq 80\%$ on most instances and obtaining $f_{\mathrm{best}}$ within $10^{-10}$ of $f^*$. Moreover, DFDS dominates PRS and IHR by its scalability, continuously attaining higher success rates and solution precision under $50\%$ lower budgets compared to PRS and IHR as dimensionality increases. 

 \section{Conclusions}
 \label{section conclusion}
In this paper, we introduce the depth-first directional search (DFDS) algorithm for the global optimization of nonconvex problems. Unlike existing methods such as hit-and-run, which primarily employ a broad-coverage search strategy, DFDS adopts a depth-first approach by performing an exhaustive line search along each sampled direction till reaching the boundary.

The core theoretical contribution of this work is the development of a geometric framework that models the probability of locating an optimizer as the surface area of a spherical cap. This framework directly leads to the establishment of the algorithm's convergence in probability and computational complexity. Specifically, we show that DFDS attains an $\epsilon-$accuracy with an expected number of function evaluations bounded by $O((\frac{2(D_0+R_\epsilon)}{\sqrt{3}R_\epsilon})^{N}\,\frac{D_0}{R_\epsilon}\,\sqrt{N}\,$\scalebox{1.25}{$\frac{1}{\epsilon}$}$)$. This analytical framework can be extended to other directional search methods as well. 

Numerical experiments on benchmark optimization problems demonstrate that DFDS achieves superior solution accuracy compared to other random search methods, such as pure random search (PRS) and improving hit-and-run (IHR). The advantage of DFDS becomes more pronounced as the instance dimension increases, validating its scalability and efficacy.

Future work will focus on several promising directions. First, integrating DFDS with efficient local search methods motivates performance enhancement and applications in high-dimensional problems. Second, extending the algorithm and its theoretical analysis to optimization over general nonconvex sets presents a challenge.

 \section*{Declarations}
\bmhead{Author Contributions} Y.-X. Zhang contributed to conceptualization, methodology, theoretical analysis, software and writing-original draft, writing-review \& editing. W.-X. Xing contributed to conceptualization, methodology, writing-review \& editing, supervision, funding acquisition. 
\bmhead{Conflict of interest}
The authors declare that they have no conflict of interest.

\bibliography{refs}
\end{document}